\documentclass[a4paper,12pt]{amsart}
\usepackage{amsfonts}
\usepackage{amssymb}
\usepackage{ifthen}
\usepackage{graphicx}
\nonstopmode \numberwithin{equation}{section}
\usepackage[usenames]{color}
\newtheorem{thm}{Theorem}[section]
\newtheorem{lem}{Lemma}[section]
\newtheorem{cor}[thm]{Corollary}
\newtheorem{prop}[thm]{Proposition}

\newtheorem{step}{Step}[section]

\theoremstyle{definition}
\newtheorem{mlem}{Main lemma}[section]
\newtheorem{assertion}{Assertion}[section]
\newtheorem{cl}{Claim}[section]
\newtheorem{ca}{Case}[section]
\newtheorem{sca}{Subcase}[section]
\newtheorem{scl}{Subclaim}[section]
\newtheorem{conj}[thm]{Conjecture}
\newtheorem{fact}{Fact}[section]
\newtheorem{defn}[thm]{Definition}
\newtheorem{op}[thm]{Open Problem}
\newtheorem{ques}[thm]{Question}
\newtheorem{rem}[thm]{Remark}
\newtheorem{exam}[thm]{Example}

\numberwithin{equation}{section}

\newcounter {own}
\def\theown {\thesection       .\arabic{own}}

\newenvironment{pf}[1][]{%
 \vskip 3mm
 \noindent
 \ifthenelse{\equal{#1}{}}%
  {{\slshape Proof. }}%
  {{\slshape #1.} }%
 }%
{\qed\bigskip}

\newcounter{alphabet}
\newcounter{tmp}
\newenvironment{Thm}[1][]{\refstepcounter{alphabet}%
\bigskip%
\noindent%
{\bf Theorem \Alph{alphabet}}%
\ifthenelse{\equal{#1}{}}{}{ (#1)}%
{\bf .} \itshape}{\vskip 8pt}

\makeatletter
\newcommand{\Ref}[1]{\@ifundefined{r@#1}{}{\setcounter{tmp}{\ref{#1}}\Alph{tmp}}}
\makeatother

\newcounter{alphabet2}

\newcommand{\ID}{{\mathbb D}}

\def\be{\begin{equation}}
\def\ee{\end{equation}}

\newcommand{\ben}{\begin{enumerate}}
\newcommand{\een}{\end{enumerate}}

\newcommand{\blem}{\begin{lem}}
\newcommand{\elem}{\end{lem}}
\newcommand{\bthm}{\begin{thm}}
\newcommand{\ethm}{\end{thm}}
\newcommand{\bcor}{\begin{cor}}
\newcommand{\ecor}{\end{cor}}
\newcommand{\beg}{\begin{exam}}
\newcommand{\eeg}{\end{exam}}
\newcommand{\begs}{\begin{examples}}
\newcommand{\eegs}{\end{examples}}
\newcommand{\bdefe}{\begin{defn}}
\newcommand{\edefe}{\end{defn}}
\newcommand{\bprob}{\begin{prob}}
\newcommand{\eprob}{\end{prob}}
\newcommand{\bques}{\begin{ques}}
\newcommand{\eques}{\end{ques}}
\newcommand{\bei}{\begin{itemize}}
\newcommand{\eei}{\end{itemize}}
\newcommand{\bcon}{\begin{conj}}
\newcommand{\econ}{\end{conj}}
\newcommand{\bop}{\begin{op}}
\newcommand{\eop}{\end{op}}

\newcommand{\bas}{\begin{assertion}}
\newcommand{\eas}{\end{assertion}}

\newcommand{\bfa}{\begin{fact}}
\newcommand{\efa}{\end{fact}}

\newcommand{\bca}{\begin{ca}}
\newcommand{\eca}{\end{ca}}

\newcommand{\bst}{\begin{step}}
\newcommand{\est}{\end{step}}

\newcommand{\bsca}{\begin{sca}}
\newcommand{\esca}{\end{sca}}

\newcommand{\bcl}{\begin{cl}}
\newcommand{\ecl}{\end{cl}}

\newcommand{\bmlem}{\begin{mlem}}
\newcommand{\emlem}{\end{mlem}}

\newcommand{\bscl}{\begin{scl}}
\newcommand{\escl}{\end{scl}}

\newcommand{\bcons}{\begin{conjs}}
\newcommand{\econs}{\end{conjs}}

\newcommand{\bprop}{\begin{prop}}
\newcommand{\eprop}{\end{prop}}

\newcommand{\br}{\begin{rem}}
\newcommand{\er}{\end{rem}}
\newcommand{\brs}{\begin{rems}}
\newcommand{\ers}{\end{rems}}
\newcommand{\bo}{\begin{obser}}
\newcommand{\eo}{\end{obser}}
\newcommand{\bos}{\begin{obsers}}
\newcommand{\eos}{\end{obsers}}
\newcommand{\bpf}{\begin{pf}}
\newcommand{\epf}{\end{pf}}
\newcommand{\ba}{\begin{array}}
\newcommand{\ea}{\end{array}}
\newcommand{\beq}{\begin{eqnarray}}
\newcommand{\beqq}{\begin{eqnarray*}}
\newcommand{\eeq}{\end{eqnarray}}
\newcommand{\eeqq}{\end{eqnarray*}}

\newcommand{\ds}{\displaystyle}

\newcounter{minutes}
\divide\time by 60
\newcounter{hours}
\multiply\time by 60 \addtocounter{minutes}{-\time}

\begin{document}

\bibliographystyle{amsplain}
\title []
{On asymptotically sharp bi-Lipschitz inequalities of quasiconformal
mappings  satisfying inhomogeneous polyharmonic equations}

\def\thefootnote{}
\footnotetext{ \texttt{\tiny File:~\jobname .tex,
          printed: \number\day-\number\month-\number\year,
          \thehours.\ifnum\theminutes<10{0}\fi\theminutes}
} \makeatletter\def\thefootnote{\@arabic\c@footnote}\makeatother

\author{Shaolin Chen}
 \address{Sh. Chen, College of Mathematics and
Statistics, Hengyang Normal University, Hengyang, Hunan 421008,
People's Republic of China.} \email{mathechen@126.com}

\author{David Kalaj}
\address{D. Kalaj, Faculty of Natural Sciences and Mathematics,
University of Montenegro, Cetinjski put b. b. 81000 Podgorica,
Montenegro. } \email{davidk@ac.me; davidkalaj@gmail.com}






\begin{abstract} For two constants $K\geq1$ and $K'\geq0$, suppose that $f$ is a $(K,K')$-quasiconformal  self-mapping of the unit
disk $\mathbb{D}$, which satisfies the following: $(1)$ the
inhomogeneous polyharmonic
  equation $\Delta^{n}f=\Delta(\Delta^{n-1}
f)=\varphi_{n}$ in $\mathbb{D}$ $(\varphi_{n}\in
\mathcal{C}(\overline{\mathbb{D}}))$, (2) the boundary conditions
$\Delta^{n-1}f=\varphi_{n-1},~\ldots,~\Delta^{1}f=\varphi_{1}$ on
$\mathbb{T}$ ($\varphi_{j}\in\mathcal{C}(\mathbb{T})$ for
$j\in\{1,\ldots,n-1\}$ and $\mathbb{T}$ denotes the unit circle),
and $(3)$ $f(0)=0$, where $n\geq2$ is an integer. The main aim of
this paper is to prove that $f$ is Lipschitz continuous, and,
further, it is bi-Lipschitz continuous when
$\|\varphi_{j}\|_{\infty}$
 are small enough for  $j\in\{1,\ldots,n\}$. Moreover, the estimates are
asymptotically sharp as $K\to 1^{+}$, $K'\to0^{+}$  and
$\|\varphi_{j}\|_{\infty}\to 0^{+}$ for  $j\in\{1,\ldots,n\}$.
\end{abstract}


\maketitle \pagestyle{myheadings} \markboth{ Sh. Chen and D. Kalaj}{
Bi-Lipschitz continuity of quasiconformal self-mappings of the unit
disk}

\section{Preliminaries and  main results }\label{csw-sec1}
Let $\mathbb{C}  \cong \mathbb{R}^{2}$ be the complex plane. For
$a\in\mathbb{C}$ and  $r>0$, let $\mathbb{D}(a,r)=\{z:\, |z-a|<r\}$, the
open disk with center $a$ and radius $r$. For convenience, we use
$\mathbb{D}_r$ to denote $\mathbb{D}(0,r)$, and $\mathbb{D}$ the
open unit disk $\mathbb{D}_1$. Let $\mathbb{T}$ be the unit circle, i.e., the
boundary $\partial\mathbb{D}$ of $\mathbb{D}$ and
$\overline{\mathbb{D}}=\mathbb{D}\cup \mathbb{T}$. Also, we denote by
$\mathcal{C}^{m}(D)$ the set of all complex-valued $m$-times
continuously differentiable functions from $D$ into $\mathbb{C}$,
where $D$ is a subset of $\mathbb{C}$  and
$m\in\{0,1,2,\ldots\}$. In particular, let
$\mathcal{C}(D):=\mathcal{C}^{0}(D)$, the set of all continuous
functions in $D$.

For a real $2\times2$ matrix $A$, we denote the matrix norm by
$\|A\|=\sup_{|z|=1}|Az|$ and the matrix function by
$\lambda(A)=\inf_{|z|=1} |Az|,$ respectively.

For $z=x+iy\in\mathbb{C}$, the
formal derivative of a complex-valued function $f=u+iv$ is given
by
$$D_{f}=\left(\begin{array}{cccc}
\ds u_{x}\;~~ u_{y}\\[2mm]
\ds v_{x}\;~~ v_{y}
\end{array}\right),
$$ where $x,~y\in\mathbb{R}$, and $u,~v$ are real-valued functions with partial derivatives.
Then,
$$\|D_{f}\|=|f_{z}|+|f_{\overline{z}}| ~\mbox{ and }~ \lambda(D_{f})=\big| |f_{z}|-|f_{\overline{z}}|\big |,
$$
where $f_{z}:=\partial f/\partial z=\frac{1}{2}\big( f_x-if_y\big)$
and $f_{\overline{z}}:=\partial f/\partial
\overline{z}=\frac{1}{2}\big(f_x+if_y\big).$ Moreover, we use
$$J_{f}:=\det D_{f} =|f_{z}|^{2}-|f_{\overline{z}}|^{2}
$$
to denote the {\it Jacobian} of $f$.

A sense-preserving homeomorphism
$$f:~\Omega_{1}\rightarrow\Omega_{2},$$ where $\Omega_{1}$ and
$\Omega_{2}$ are subdomains of $\mathbb{C}$, is said to be a {\it
$(K,K')$-quasiconformal} mapping if $f$ is absolutely continuous on
lines in $\Omega_{1}$, and there are constants $K\geq1$ and
$K'\geq0$ such that $$\|D_{f}(z)\|^{2}\leq
KJ_{f}(z)+K',~z\in\Omega_{1}.$$ In particular, if $K'=0$, then $f$
is a $K$-quasiconformal mapping (cf. \cite{FS,Ni}).

Given a subset $\Omega$ of $\mathbb{C}$, a function
$\psi:~\Omega\rightarrow\mathbb{C}$ is said to be {\it bi-Lipschitz}
if there is a  constant $L\geq1$ such that for all
 $z_{1},z_{2}\in\Omega$,

\be\label{eq-L}\frac{1}{L}|z_{1}-z_{2}|\leq|\psi(z_{1})-\psi(z_{2})|\leq
L|z_{1}-z_{2}|.\ee In particular,  $\psi$ is called {\it Lipschitz}
if the inequality on the right  of (\ref{eq-L}) holds, and $\psi$ is
said to be {\it co-Lipschitz} if it satisfies the inequality on the
left  of (\ref{eq-L}). It is clear that any sense-preserving
bi-Lipschitz mapping is quasiconformal mapping (see Chapter 14.78 in
\cite{HKM}). But quasiconformal mappings are not necessarily
bi-Lipschitz, not even Lipschitz (see \cite{KS}).










In \cite{Pav2}, Pavlov\'c discussed the bi-Lipschitz characteristic
of a harmonic homeomorphism of $\mathbb{D}$ onto itself.  Later,
Partyka and Sakan \cite{PS} gave explicit estimates of bi-Lipschitz
constants for a harmonic $K$-quasiconformal mapping $f$ of
$\mathbb{D}$ onto itself. Under the additional assumption $f(0)=0$,
the estimates are asymptotically sharp as $K\rightarrow1^{+}$, so
$f$ behaves almost like a rotation for sufficiently small $K$.
Recently, the  bi-Lipschitz characteristics of harmonic
quasiconformal mappings have been attracted much attention (see
\cite{K4, K7, KMa, MV1, PS, Pav2}). The Lipschitz continuity of $(K,
K')$-quasiconformal harmonic mappings  has also been investigated in
\cite{CLSW,K8,ZTY}. On the discussion of the related topic, we refer
to \cite{Ada,Ale,H-S, K3, K9, Man,Ma, Pav2} and the related references therein.
On the study of the Lipschitz characteristic of quasiconformal
mappings satisfying certain elliptic PDEs, see
\cite{CW,K6,K-M,K2,KS}. In particular, let us recall the following
two results.

\begin{Thm} {\rm (\cite[Theorem 1.1]{CW})}\label{Thm-1.1}
Let $g\in\mathcal{C}(\overline{\mathbb{D}})$,
$\varphi\in\mathcal{C}(\mathbb{T})$, and let $K\geq1$ be a constant.
Suppose that $f$ is a $K$-quasiconformal self-mapping of
$\mathbb{D}$ satisfying the bi-harmonic equation $\Delta(\Delta
f)=g$ with $\Delta f=\varphi$ in $\mathbb{T}$ and $f(0)=0$. Then,
there are nonnegative constants $\mathcal{M}_{j}(K)$ and
$\mathcal{N}_{j}(K,\varphi,g)$ {\rm($j\in\{1,2\}$)} with
$\lim_{K\rightarrow1^{+}}\mathcal{M}_{j}(K)=1$ and
$\lim_{\|\varphi\|_{\infty}\rightarrow0^{+},\|g\|_{\infty}\rightarrow0^{+}}\mathcal{N}_{j}(K,\varphi,g)=0$
such that for all $z_{1}$ and $z_{2}$ in $\ID$,
$$
\big(\mathcal{M}_{1}(K)-\mathcal{N}_{1}(K,\varphi,g)\big)|z_{1}-z_{2}|\leq|f(z_{1})-f(z_{2})|\leq
\big(\mathcal{M}_{2}(K)+\mathcal{N}_{2}(K,\varphi,g)\big)|z_{1}-z_{2}|,$$
where
$\|\varphi\|_{\infty}=\sup_{\zeta\in\mathbb{T}}|\varphi(\zeta)|$ and
$\|g\|_{\infty}=\sup_{z\in\mathbb{D}}|g(z)|$.
\end{Thm}

\begin{Thm} {\rm (\cite[Corollary 3.1]{CW})}\label{Thm-1.2}
Under the assumptions of Theorem \Ref{Thm-1.1}, if, further,
$\|g\|_{\infty}\leq a_1(K)$ and $\|\varphi\|_{\infty}\leq a_2(K)$,
then $f$ is co-Lipschitz continuous, and so, it is bi-Lipschitz
continuous, where $a_1(K)=\frac{60}{(25+61K^{2})46^{2(K-1)}}$ and
$a_2(K)=\frac{25}{(38+101K^{2})46^{2(K-1)}}.$
\end{Thm}

We remark that  the bi-harmonic quasiconformal mappings between
smooth domains are not necessarily Lipschitz (see the  Example
\ref{CKW-1}).
\begin{exam}\label{CKW-1}
The mapping $f(z)=z \log(|z|^2)$ is {\it bi-harmonic} (i.e.,
$\Delta(\Delta f)=0$) in $\mathbb{D}_{e^{-2}}\setminus\{0\}$ and
quasiconformal in $\mathbb{D}_{e^{-2}}$
 but is not Lipschitz in any neighborhood of $z=0$. The mapping $f$ is not
bi-harmonic  in $0$, since $\Delta f(z)=1/\bar z$.  The mapping $f(z)=z\log(|z|^2)$ is bi-harmonic in $\mathbb{D}(ri, r)$ 
for small enough positive number $r$, and maps $\mathbb{D}(ri, r)$
onto a convex Jordan domain $\Omega$ with $\mathcal{C}^2$ boundary.
Thus the bi-harmonic mapping $h(z)=f(r(z+ i))$ maps $\mathbb{D}$
quasiconformally onto the Jordan domain $\Omega$ with
$\partial\Omega\in \mathcal{C}^2$, but it is not Lipschitz. So the
Lipschitz continuity fails if we drop the condition  $\Delta f$ is
continuous up to the boundary. To prove that $\partial\Omega\in
\mathcal{C}^2$ we observe first that $\partial\Omega$ is
rectifiable. Namely by direct computation, we have
\[\begin{split}|\partial
\Omega|&=\int_0^{2\pi}\left|\frac{\partial}{\partial t}
h(e^{it})\right|dt\\&=\int_0^2\pi r \sqrt{1+\left(1+\log\left[2 r^2
(1+\sin t)\right]\right)^2-2 \sin t-2\sin t \log\left[2 r^2 (1+\sin
t)\right]}dt\\&<\infty \end{split}\] and
\begin{equation}\label{equ}\frac{\frac{\partial}{\partial t} h(e^{it})}{\left|\frac{\partial}{\partial t} h(e^{it})\right|}=\frac{{1+i e^{i t} (1+\log
[r (1+\sin t )])}}{|{1+i e^{i t} (1+\log [r (1+\sin t )])}|}.
\end{equation}
Since $$\frac{\frac{\partial}{\partial t}
h(e^{it})}{\left|\frac{\partial}{\partial t} h(e^{it})\right|}
=e^{i\varphi(s(t))},$$ where $$s(t)=\int_0^t
\left|\frac{\partial}{\partial \tau} h(e^{i\tau})\right|d\tau$$ is
the natural parameter, and since the limit of left-hand side in
\eqref{equ} tends to $-1$, it follows that $\varphi$ is continuous
in $s=0$, and therefore the function $s\to h(t(s))$ is
$\mathcal{C}^1$. To show that the curve is $\mathcal{C}^2$, we find
the curvature of $\partial\Omega$ at $0$. Namely if $x(t)=\mbox{Re}
(h(e^{it}))$ and $y(t)=\mbox{Im} (h(e^{it}))$, then the curvature
$$\kappa(t)=\frac{|\ddot x \dot y-\ddot y\dot x|}{(\dot x^2+\dot y^2)^{3/2}}.$$ Then it can be proved that $\lim_{t\to 0}\kappa(t)=0$.
 Thus $\kappa$ is continuous in $\partial\Omega$ which means that the curve is $\mathcal{C}^2$.
\end{exam}

The main aim of this paper is to improve and generalize Theorems
\Ref{Thm-1.1} and \Ref{Thm-1.2}. In order to state our main results,
we need to recall some basic definitions and some results which
motivate the present work.

For $z, \zeta\in\mathbb{D}$ with $z\neq \zeta$, let
$$G(z,\zeta)=\frac{1}{2\pi}\log\left|\frac{1-z\overline{\zeta}}{z-\zeta}\right|~\mbox{ and}
~P(z,e^{it})=\frac{1}{2\pi}\frac{1-|z|^{2}}{|1-ze^{-it}|^{2}}$$ be
the {\it
 Green function} and   {\it (harmonic) Poisson kernel},
respectively, where $t\in[0,2\pi].$

Let $\varphi_{n}\in \mathcal{C}(\overline{\mathbb{D}})$ and
$f\in\mathcal{C}^{2n}(\mathbb{D})$, where $n\geq2$ is an integer. Of
particular interest for our investigation is the following {\it
inhomogeneous polyharmonic equation (or n-harmonic equation)}:

\be\label{eq-ch-1} \Delta^{n}f=\Delta(\Delta^{n-1}
f)=\varphi_{n}~\mbox{in}~\mathbb{D}\ee with the following associated
{\it Dirichlet boundary value condition}:

\be\label{eq-ch-2}
\Delta^{n-1}f=\varphi_{n-1},~\ldots,~\Delta^{1}f=\varphi_{1},~\Delta^{0}f=\varphi_{0}~\mbox{on}~\mathbb{T},\ee
where $\Delta^{0}f:=f$,
$$\Delta^{1} f:=\Delta f=\frac{\partial^{2}f}{\partial
x^{2}}+\frac{\partial^{2}f}{\partial y^{2}}=4f_{z \overline{z}}$$
stands for the {\it Laplacian} of $f$, and  $\varphi_{k}\in
\mathcal{C}(\mathbb{T})$ for $k\in\{0,1,\ldots,n-1\}$. Here
$\Delta^{j-1}f=\varphi_{j-1}$ on $\mathbb{T}$ means that
$$\lim_{r\rightarrow1}\Delta^{j-1}f(r\zeta)=\varphi_{j-1}(\zeta)$$
for all $\zeta\in\mathbb{T}$, where $j\in\{1,2,\ldots,n\}$.

By the iterated {\it Poly-Cauchy integral operators} (cf.
\cite{Be}), we see that all solutions to the equation
(\ref{eq-ch-1}) satisfying (\ref{eq-ch-2}) are given by

\be\label{eq-ch-3.0}f(z)=P[\varphi_{0}](z)+\sum_{k=1}^{n}(-1)^{k}G_{k}[\varphi_{k}](z),~z\in\mathbb{D},\ee
where
$$P[\varphi_{0}](z)=\int_{0}^{2\pi}P(z,e^{it})\varphi_{0}(e^{it})dt,$$

\beq\label{eq-ch-3.1} G_{k}[\varphi_{k}](z)&=&
\int_{\mathbb{D}}\cdots\int_{\mathbb{D}}G(z,\xi_{1})\cdots
G(\xi_{k-1},\xi_{k})\\ \nonumber
&&\times\left(\int_{0}^{2\pi}P(\xi_{k},e^{it})\varphi_{k}(e^{it})dt\right)d\sigma(\xi_{k})\cdots
d\sigma(\xi_{1})\eeq for $k\in\{1,\ldots,n-1\}$, and

\beq\label{eq-ch-3.2} G_{n}[\varphi_{n}](z)&=&
\int_{\mathbb{D}}\cdots\int_{\mathbb{D}}G(z,\zeta_{1})\cdots G(\zeta_{n-2},\zeta_{n-1})\\
\nonumber
&&\times\left(\int_{\mathbb{D}}G(\zeta_{n-1},\zeta_{n})\varphi_{n}(\zeta_{n})d\sigma(\zeta_{n})\right)d\sigma(\zeta_{n-1})\cdots
d\sigma(\zeta_{1}).\eeq Here $d\sigma$ is the Lebesgue area measure
in $\mathbb{D}$. The behavior of  solutions to the  polyharmonic
equations with the different boundary value conditions has attracted
much attention of many authors (cf.
\cite{GGS,GR-1997,HK,Lai,May-M,OL}).


The first aim of this paper is to investigate the  asymptotically
sharp bi-Lipschitz inequalities of $(K,K')$-quasiconformal
self-mapping of $\mathbb{D}$ satisfying the inhomogeneous
polyharmonic equation \eqref{eq-ch-1} with the boundary condition
\eqref{eq-ch-2}. It is read as follows.

\begin{thm}\label{thm-1.12} Let $\varphi_{0}$ be a sense-preserving
homeomorphism of $\mathbb{T}$ onto itself. For $n\geq2$ and
$k\in\{1,\ldots,n-1\}$, let
$\varphi_{n}\in\mathcal{C}(\overline{\mathbb{D}})$ and
$\varphi_{k}\in\mathcal{C}(\mathbb{T})$, and let $K\geq1$ and
$K'\geq0$ be  constants. Suppose that $f$ is a
$(K,K')$-quasiconformal self-mapping of $\overline{\mathbb{D}}$
satisfying the inhomogeneous polyharmonic equation \eqref{eq-ch-1}
with the Dirichlet boundary value condition \eqref{eq-ch-2}.

\begin{enumerate}
\item[{\rm $(a)$}] If
$|P[\varphi_{0}](0)|+\sqrt{K'}+2K\left(\frac{1}{3}\|\varphi_{1}\|_{\infty}+\frac{1}{15}\sum_{k=2}^{n}\left(\frac{3}{16}\right)^{k-2}\|\varphi_{k}\|_{\infty}\right)
<\frac{2}{\pi},$ then $f$ is bi-Lipschitz continuous in
$\mathbb{D}$.
\item[{\rm $(b)$}] If $P[\varphi_{0}](0)=0$ and $\sqrt{K'}+2K\left(\frac{1}{3}\|\varphi_{1}\|_{\infty}+\frac{1}{15}
\sum_{k=2}^{n}\left(\frac{3}{16}\right)^{k-2}\|\varphi_{k}\|_{\infty}\right)
<\frac{2}{\pi},$ then there are nonnegative constants $M_{j}(K,K')$
and $N_{j}(K,\varphi_{1},\cdots,\varphi_{n})$ {\rm($j\in\{1,2\}$)}
with
$$\lim_{K\rightarrow1^{+},K'\rightarrow0^{+}}M_{j}(K,K')=1,~\lim_{\|\varphi_{1}\|_{\infty}\rightarrow0^{+},\cdots,\|\varphi_{n}\|_{\infty}\rightarrow0^{+}}N_{j}(K,\varphi_{1},\cdots,\varphi_{n})=0$$
and $M_{2}(K,K')-N_{2}(K,\varphi_{1},\cdots,\varphi_{n})>0$ such
that for all $z_{1}, z_{2}\in\mathbb{D}$,
\begin{eqnarray*}
\big(M_{2}(K,K')&-&N_{2}(K,\varphi_{1},\cdots,\varphi_{n})\big)|z_{1}-z_{2}|\leq|f(z_{1})-f(z_{2})|\\&\leq&
\big(M_{1}(K,K')+N_{1}(K,\varphi_{1},\cdots,\varphi_{n})\big)|z_{1}-z_{2}|,
\end{eqnarray*}
where
$\|\varphi_{n}\|_{\infty}=\sup_{z\in\mathbb{D}}|\varphi_{n}(z)|$ and
$\|\varphi_{k}\|_{\infty}=\sup_{\zeta\in\mathbb{T}}|\varphi_{k}(\zeta)|$
for $k\in\{1,\ldots,n-1\}$.
\end{enumerate}

\end{thm}



In particular, if $K'=0$, then we have the following better
estimate.


\begin{thm} \label{thm-1.1}
Let $\varphi_{n}\in\mathcal{C}(\overline{\mathbb{D}})$ and
$\varphi_{k}\in\mathcal{C}(\mathbb{T})$, and let $K\geq1$ be a
constant, where $n\geq2$ and $k\in\{1,\ldots,n-1\}$. Suppose that
$f$ is a $K$-quasiconformal self-mapping of $\mathbb{D}$ satisfying
the inhomogeneous polyharmonic equation \eqref{eq-ch-1} with
$\Delta^{n-1}f=\varphi_{n-1},~\ldots,~\Delta^{1}f=\varphi_{1}$ on
$\mathbb{T}$ and $f(0)=0$. Then, there are nonnegative constants
$M_{j}(K)$ and $N_{j}(K,\varphi_{1},\cdots,\varphi_{n})$
{\rm($j\in\{3,4\}$)} with
$$\lim_{K\rightarrow1^{+}}M_{j}(K)=1~\mbox{and}~\lim_{\|\varphi_{1}\|_{\infty}\rightarrow0^{+},\cdots,\|\varphi_{n}\|_{\infty}\rightarrow0^{+}}N_{j}(K,\varphi_{1},\cdots,\varphi_{n})=0$$
such that for all $z_{1},z_{2}\in\mathbb{D}$,
\begin{eqnarray*}
\big(M_{4}(K)-N_{4}(K,\varphi_{1},\cdots,\varphi_{n})\big)|z_{1}-z_{2}|&\leq&|f(z_{1})-f(z_{2})|\\&\leq&
\big(M_{3}(K)+N_{3}(K,\varphi_{1},\cdots,\varphi_{n})\big)|z_{1}-z_{2}|.
\end{eqnarray*}
\end{thm}



The following is the so-called Mori's Theorem (cf. \cite{CL, K2,
Mo}). We refer to \cite{FV,MN} for some  analogical results of
Theorem \Ref{Mori} in the higher dimensional case.

\begin{Thm}\label{Mori}
Suppose that $f$ is a $K$-quasiconformal self-mapping of
$\mathbb{D}$ with $f(0)=0$. Then, there exists a constant $Q(K)$,
satisfying the condition $Q(K)\rightarrow1$ as $K\rightarrow1$, such
that
$$|f(z_{2})-f(z_{1})|\leq Q(K)|z_{2}-z_{1}|^{\frac{1}{K}},$$
where the notation $Q(K)$ means that the constant $Q$ depends only
on $K$.
\end{Thm}

We remark that in \cite{Qiu} it is proved

\be\label{eq-qiu}1\leq
Q(K)\leq16^{1-\frac{1}{K}}\min\left\{\bigg(\frac{23}{8}\bigg)^{1-\frac{1}{K}},~\big(1+2^{3-2K}\big)^{\frac{1}{K}}\right\}.\ee

As a direct consequence of Claim \ref{claim-3.8} in the proof of
Theorem \ref{thm-1.1}, we have the following result.

\bcor\label{tue-1} Under the assumptions of Theorem \ref{thm-1.1},
if, further,

\begin{eqnarray*}
\frac{(Q(K))^{-2K}K^{-2}}{2\pi}\int_{0}^{2\pi}|e^{it}-e^{i\theta}|^{2K-2}dt&>&
\left(\frac{7}{6}+\frac{1}{2K^{2}}\right)\|\varphi_{1}\|_{\infty}
\\&&+\sum_{j=2}^{n}\left(\frac{47}{240}+\frac{1}{16K^{2}}\right)\|\varphi_{j}\|_{\infty}\left(\frac{3}{16}\right)^{j-2},
\end{eqnarray*}
then $f$ is co-Lipschitz continuous, and so, it is bi-Lipschitz
continuous, where $Q(K)$ is the same as in Theorem \Ref{Mori}. \ecor

By (\ref{eq-qiu}) and \cite[Formula 3.27]{K2}, we see that
\begin{eqnarray*}
K^{-2}(Q(K))^{-2K}\frac{1}{2\pi}\int_{0}^{2\pi}|e^{it}-e^{i\theta}|^{2K-2}dt&=&\frac{2^{2K-2}\Gamma\big(K-\frac{1}{2}\big)}{\sqrt{\pi}K^{2}(K-1)\Gamma(K-1)(Q(K))^{2K}}\\
&\geq&\frac{1}{K^{2}(Q(K))^{2K}}\geq\frac{1}{K^{2}46^{2K-2}},
\end{eqnarray*}
which gives the following result, where  $\Gamma$ is the Gamma
function.

\bcor\label{tue-2} Under the assumptions of Theorem \ref{thm-1.1},
if, further,
$$
\frac{1}{K^{2}46^{2K-2}}>
\left(\frac{7}{6}+\frac{1}{2K^{2}}\right)\|\varphi_{1}\|_{\infty}
+\sum_{j=2}^{n}\left(\frac{47}{240}+\frac{1}{16K^{2}}\right)\|\varphi_{j}\|_{\infty}\left(\frac{3}{16}\right)^{j-2},
$$
then $f$ is co-Lipschitz continuous, and so, it is bi-Lipschitz
continuous. \ecor

We remark that if $n=2$, then Corollary \ref{tue-2} is an
improvement of Theorem \Ref{Thm-1.2}.

By the discussions in Step 3 of the proof of Theorem \ref{thm-1.1}
in Section \ref{sec-4} or by Corollary \ref{tue-1}, we see that the
co-Lipschitz continuity coefficient
$$M_{4}(K)-N_{4}(K,\varphi_{1},\cdots,\varphi_{n})$$ is positive for
small enough norms  $\|\varphi_{k}\|_{\infty}$, where
$k\in\{1,\ldots,n\}$. The following example (Example \ref{wed-5})
shows that the condition for $f$ to be co-Lipschitz continuous cannot
be replaced by the one that $\varphi_{k}$  are arbitrary, where
$k\in\{1,\ldots,n\}$.

\begin{exam}\label{wed-5}
 For $z\in\overline{\mathbb{D}}$, let
 $$\varphi_{n}(z)=\beta\left(\prod_{j=1}^{n}(\tau-2j+4)\right)\left(\prod_{j=1}^{n}(\tau-2j+2)\right)z^{\frac{\tau}{2}-n+1}\overline{z}^{\frac{\tau}{2}-n},$$
  where
$\tau>2n-1$  and $\beta$ are constants with $n\geq2$ and
$|\beta|=1$.  Suppose that $f$ satisfies the following polyharmonic equation

\be\label{ckw-p1}\Delta^{n} f(z)=\Delta(\Delta^{n-1}
f(z))=\varphi_{n}(z),~z\in\mathbb{D},\ee
 with the following associated Dirichlet boundary value condition:
 $$\Delta^{k}
f(\xi)=\varphi_{k}(\xi)~\mbox{and}~f(\xi)=\varphi_{0}(\xi),~\xi\in\mathbb{T},$$
where $\varphi_{0}(\xi)=\beta\xi$, and for $k\in\{1,\ldots,n-1\}$, $$\varphi_{k}(\xi)=\beta\left(\prod_{j=1}^{k}(\tau-2j+4)\right)\left(\prod_{j=1}^{k}(\tau-2j+2)\right)\xi.$$

It follows from (\ref{eq-ch-3.0}) that
$$f(z)=\beta|z|^{\tau}z,~z\in\overline{\mathbb{D}},$$ is the solution to (\ref{ckw-p1}).
 Obviously,  $f$ is a  $K$-quasiconformal self-mapping of $\mathbb{D}$ with $f(0)=0$
and $K=1+\tau$. Furthermore,
$$\|\varphi_{n}\|_{\infty}=\left(\prod_{j=1}^{n}(\tau-2j+4)\right)\left(\prod_{j=1}^{n}(\tau-2j+2)\right),$$
and for $k\in\{1,\ldots,n-1\}$,
$$\|\varphi_{k}\|_{\infty}=\left(\prod_{j=1}^{k}(\tau-2j+4)\right)\left(\prod_{j=1}^{k}(\tau-2j+2)\right) \;\;\mbox{and}\;\; \|\varphi_{0}\|_{\infty}=1.$$
However, $f$ is not co-Lipschitz continuous because
$$\lambda(D_{f}(0))=|f_{z}(0)|-|f_{\overline{z}}(0)|=0.$$
\end{exam}


By applying Corollary \ref{tue-2}, we illustrate the possibility of
$f$ to be bi-Lipschitz continuous by the following example.

\begin{exam}\label{wed-6}
Suppose that $f$ satisfies the following bi-harmonic equation

\be\label{ckw-p2}\Delta(\Delta f(z))=-\frac{16}{15},~z\in\mathbb{D},\ee
with the following associated Dirichlet boundary value condition:
$$\Delta f(\xi)=-\frac{1}{5}~\mbox{and}~f(\xi)=\xi,~\xi\in\mathbb{T}.$$

By (\ref{eq-ch-3.0}), we see that

$$f(z)=z+\frac{1}{60}(|z|^{2}-|z|^{4}),~z\in\overline{\mathbb{D}},$$
is the solution to (\ref{ckw-p2}).
It is not difficult to know that $f$ is a
$K$-quasiconformal self-mapping of $\mathbb{D}$ with


$$
K=\max_{z\in\overline{\mathbb{D}}}\left\{\frac{|1+\overline{z}(1-2|z|^{2})M|+|Mz(1-2|z|^{2})|}{|1+\overline{z}(1-2|z|^{2})M|-|Mz(1-2|z|^{2})|}\right\}
=\frac{30}{29},$$ where $M=\frac{1}{60}$. Since elementary
computations lead to
$$\frac{1}{K^{2}46^{2(K-1)}}=\frac{29^{2}}{30^{2}46^{\frac{2}{29}}}>0.717,~\left(\frac{7}{6}+\frac{1}{2K^{2}}\right)\|\varphi_{1}\|_{\infty}<0.326$$
and
$$\left(\frac{47}{240}+\frac{1}{16K^{2}}\right)\|\varphi_{2}\|_{\infty}<0.271,$$
we see that
$$
\frac{1}{K^{2}46^{2K-2}}>
\left(\frac{7}{6}+\frac{1}{2K^{2}}\right)\|\varphi_{1}\|_{\infty}
+\left(\frac{47}{240}+\frac{1}{16K^{2}}\right)\|\varphi_{2}\|_{\infty},
$$ where $\|\varphi_{1}\|_{\infty}=\frac{1}{5}$ and $\|\varphi_{2}\|_{\infty}=\frac{16}{15}.$
Now, it follows from Corollary \ref{tue-2} that $f$ is co-Lipschitz
continuous, and so, it is bi-Lipschitz continuous.
\end{exam}



We recall that the (periodic) {\it Hilbert transformation}  of a
$2\pi-$periodic function $\Psi\in L^{1}(\mathbb{T})$ is defined by
$$ H(\Psi)(\theta)=-\frac{1}{\pi}\int_{0}^{\pi}\frac{\Psi(\theta+t)-\Psi(\theta+t)}{2\tan(t/2)}dt.$$
It is well known that the Lipschitz continuity of  $\varphi$ in
$\mathbb{T}$ is not enough to guarantee that its harmonic extension
$P[\varphi]$ is also Lipschitz continuous. In fact, $P[\varphi]$ is
Lipschitz continuous if and only if the Hilbert transform of
$d\varphi(e^{i\theta})/d\theta\in L^{\infty}(\mathbb{T})$ (cf.
\cite{CLW,Z}). The last aim of this paper is to investigate the
Lipschitz continuity of   solutions to the inhomogeneous
polyharmonic equation (\ref{eq-ch-1}) satisfying some certain
boundary conditions.

\begin{prop}\label{thm-02}
For $n\geq2$ and $k\in\{1,\ldots,n-1\}$, let
$\varphi_{n}\in\mathcal{C}(\overline{\mathbb{D}})$ and
$\varphi_{k}\in\mathcal{C}(\mathbb{T})$, and let
$\varphi_{0}\in\mathcal{C}(\mathbb{T})$ be differentiable. Suppose
that $f$ is a solution to the inhomogeneous polyharmonic equation
\eqref{eq-ch-1} satisfying
$\Delta^{n-1}f=\varphi_{n-1},~\ldots,~\Delta^{1}f=\varphi_{1},~\Delta^{0}f=\varphi_{0}$
on $\mathbb{T}$. Then $f$ is Lipschitz continuous in $\mathbb{D}$ if
and only if the Hilbert transform of
$d\varphi_{0}(e^{i\theta})/d\theta\in L^{\infty}(\mathbb{T})$.
\end{prop}

 We will prove several auxiliary
results in  Section \ref{sec-2}. The proof of Theorem \ref{thm-1.12}
will be presented in Section \ref{sec-3}. Theorem \ref{thm-1.1} and
Proposition \ref{thm-02} will be proved in Sections \ref{sec-4} and
\ref{sec-5}, respectively.

\section{Some auxiliary results}\label{sec-2}
In this section, we shall prove several lemmas which will be used later on. The first lemma is as follows.

\begin{lem}\label{lem-0.1} Let $G$ be the Green function. Then, for $z\in\mathbb{D}$,
\be\label{kk-1}\int_{\mathbb{D}}|G(z,\zeta)|d\sigma(\zeta)=\frac{1-|z|^{2}}{4}\ee
and
\be\label{kk-2}\int_{\mathbb{D}}(1-|\zeta|^{2})|G(z,\zeta)|d\sigma(\zeta)=\frac{(1-|z|^{2})(3-|z|^{2})}{16}\leq\frac{3(1-|z|^{2})}{16}.\ee
\end{lem}

\begin{Thm}{\rm (cf. \cite{LP})}\label{Thm-1}
For  $z\in\mathbb{D}$, we have
$$\frac{1}{2\pi}\int_{0}^{2\pi}\frac{d\theta}{|1-ze^{i\theta}|^{2\alpha}}=
\sum_{k=0}^{\infty}\left(\frac{\Gamma(k+\alpha)}{k!\Gamma(\alpha)}\right)^{2}|z|^{2k},$$
where $\alpha>0$.
\end{Thm}

\subsection*{Proof of Lemma \ref{lem-0.1}} We first prove
(\ref{kk-1}). Let \be\label{kk-4}
w=\frac{z-\zeta}{1-\overline{z}\zeta}=re^{it},\ee where $r\in[0,1)$
and $t\in[0,2\pi]$.
 Since Theorem \Ref{Thm-1} implies
\be\label{kk-5}\frac{1}{2\pi}\int_{0}^{2\pi}\frac{dt}{|1-\overline{z}re^{it}|^{4}}=\sum_{j=0}^{\infty}(j+1)^{2}|z|^{2j}r^{2j},\ee
by (\ref{kk-4}),
 we
obtain
\beqq\label{yy-1}\int_{\mathbb{D}}|G(z,\zeta)|d\sigma(\zeta)=\frac{1}{2\pi}\int_{\mathbb{D}}\left(\log\frac{1}{|w|}\right)
\frac{(1-|z|^{2})^{2}}{|1-\overline{z}w|^{4}}d\sigma(w)=\frac{(1-|z|^{2})}{4}.\eeqq

 Now we show  (\ref{kk-2}). For $z\in\mathbb{D}$,
let
$$I_{1}(z)=\int_{\mathbb{D}}(1-|\zeta|^{2})|G(z,\zeta)|d\sigma(\zeta).$$


By Theorem \Ref{Thm-1}, we have

$$\frac{1}{2\pi}\int_{0}^{2\pi}\frac{dt}{|1-\overline{z}re^{it}|^{6}}=\sum_{j=0}^{\infty}\frac{(j+1)^{2}(j+2)^{2}}{4}|z|^{2j}r^{2j},$$
which, together with
$$1-|\zeta|^{2}=\frac{(1-|w|^{2})(1-|z|^{2})}{|1-\overline{z}w|^{2}},$$
implies that



\begin{eqnarray*} I_{1}(z)&=&
\frac{1}{2\pi}\int_{\mathbb{D}}\frac{(1-|w|^{2})(1-|z|^{2})^{3}}{|1-\overline{z}w|^{6}}\log\frac{1}{|w|}d\sigma(w)\\
\nonumber
&=&(1-|z|^{2})^{3}\int_{0}^{1}(1-r^{2})r\log\frac{1}{r}\left(\frac{1}{2\pi}\int_{0}^{2\pi}\frac{dt}{|1-\overline{z}re^{it}|^{6}}\right)dr\\
\nonumber &=&\frac{(1-|z|^{2})(3-|z|^{2})}{16}.
 \end{eqnarray*}
The proof of this lemma is completed.\qed
\medskip





\begin{lem}\label{lem-es-5} For $z\in\mathbb{D}$,
\be\label{bst-5}\frac{1}{2\pi}\int_{\mathbb{D}}\frac{(1-|\varsigma|^{2})^{2}}{|1-z\overline{\varsigma}||z-\varsigma|}
d\sigma(\varsigma)\leq\frac{4(2-|z|^{2})}{15}.\ee In particular, the
inequality {\rm(\ref{bst-5})} is sharp at $z=0$.
\end{lem}

\bpf Let
$$\eta=\frac{z-\varsigma}{1-\overline{z}\varsigma}=\rho e^{i\theta}.$$
Then

\be\label{eq-y12}1-\overline{z}\varsigma=\frac{1-|z|^{2}}{1-\overline{z}\eta}~\mbox{and}~1-|\varsigma|^{2}=\frac{(1-|\eta|^{2})(1-|z|^{2})}{|1-\overline{z}\eta|^{2}}.
\ee By Theorem \Ref{Thm-1}, we obtain

\beq\label{eqs7}
\frac{1}{2\pi}\int_{\mathbb{D}}\frac{(1-|\eta|^{2})^{2}}{|\eta||1-\overline{z}\eta|^{6}}d\sigma(\eta)&=&\int_{0}^{1}(1-\rho^{2})^{2}\left(\frac{1}{2\pi}\int_{0}^{2\pi}\frac{d\theta}{|1-\overline{z}\rho
e^{i\theta}|^{6}}\right)d\rho\\ \nonumber
&=&\sum_{j=0}^{\infty}\frac{(j+1)^{2}(j+2)^{2}}{4}|z|^{2j}\int_{0}^{1}(1-\rho^{2})^{2}\rho^{2j}d\rho\\
\nonumber
&=&\sum_{j=0}^{\infty}\frac{(j+1)^{2}(j+2)^{2}}{\big(2j+1\big)\big(2j+3\big)\big(j+\frac{5}{2}\big)}|z|^{2j}.
\eeq By computation,  we have
\be\label{eq-y13}\frac{(j+1)^{2}(j+2)^{2}}{\big(2j+1\big)\big(2j+3\big)\big(j+\frac{5}{2}\big)}\leq\frac{4}{15}(j+2).\ee

It follows from   (\ref{eq-y12}),  (\ref{eqs7}) and (\ref{eq-y13})
that
\begin{eqnarray*}
\frac{1}{2\pi}\int_{\mathbb{D}}\frac{(1-|\varsigma|^{2})^{2}}{|1-z\overline{\varsigma}||z-\varsigma|}
d\sigma(\varsigma)&=&\frac{(1-|z|^{2})^{2}}{2\pi}\int_{\mathbb{D}}\frac{(1-|\eta|^{2})^{2}}{|\eta||1-\overline{z}\eta|^{6}}d\sigma(\eta)\\
&=&(1-|z|^{2})^{2}\sum_{j=0}^{\infty}\frac{(j+1)^{2}(j+2)^{2}}{\big(2j+1\big)\big(2j+3\big)\big(j+\frac{5}{2}\big)}|z|^{2j}\\
&\leq&\frac{4(1-|z|^{2})^{2}}{15}\sum_{j=0}^{\infty}(j+2)|z|^{2j}\\
&=&\frac{4(2-|z|^{2})}{15}.
\end{eqnarray*} The proof of this lemma is finished.
\epf

\begin{lem}\label{lem-y1} Let $P$ be the Poisson kernel and $\theta\in[0,2\pi]$. Then
$$\int_{\mathbb{D}}P(\zeta,e^{i\theta})(1-|\zeta|^{2})d\sigma(\zeta)=\frac{1}{4}.$$

\end{lem}
\bpf Let $\zeta=\varrho e^{it}$. By Theorem \Ref{Thm-1}, we have
$$\frac{1}{2\pi}\int_{0}^{2\pi}\frac{dt}{|1-\varrho e^{it}e^{-i\theta}|^{2}}=\frac{1}{1-\varrho^{2}},$$
which gives that
\begin{eqnarray*}
\int_{\mathbb{D}}P(\zeta,e^{i\theta})(1-|\zeta|^{2})d\sigma(\zeta)&=&\int_{0}^{1}(1-\varrho^{2})^{2}\varrho\left[\frac{1}{2\pi}\int_{0}^{2\pi}\frac{dt}{|1-\varrho
e^{it}e^{-i\theta}|^{2}}\right]d\varrho\\
&=&\frac{1}{4}.
\end{eqnarray*}
\epf

\begin{lem}\label{lem-1}
Suppose that  $\varphi_{k}\in \mathcal{C}(\mathbb{T})$  and
$G_{k}[\varphi_{k}]$ are defined in \eqref{eq-ch-3.1}, where
$k\in\{1,\ldots,n-1\}$ and $n\geq2$. Then, the following statements
hold:

\noindent $(1)$ For $z\in\mathbb{D}$,

$$\max\Bigg\{\left|\frac{\partial}{\partial
z}G_{k}[\varphi_{k}](z)\right|,\;\;\left|\frac{\partial}{\partial
\overline{z}}G_{k}[\varphi_{k}](z)\right|\Bigg\}\leq\nu_{k}(z),$$
where

$$\nu_{k}(z)=\begin{cases}
\displaystyle \frac{1}{3}\|\varphi_{1}\|_{\infty},
&\mbox{ if }\, k=1,\\
\displaystyle
\|\varphi_{k}\|_{\infty}\left(\frac{3}{16}\right)^{k-2}\frac{(2-|z|^{2})}{30},
&\mbox{if}\, 2\leq k\leq n-1.
\end{cases}$$

\noindent $(2)$ Both $\frac{\partial}{\partial z}G_{k}[\varphi_{k}]$
and $\frac{\partial}{\partial \overline{z}}G_{k}[\varphi_{k}]$ have
continuous extensions to the boundary, and further, for
$\theta\in[0,2\pi]$,
$$\max\Bigg\{\left|\frac{\partial}{\partial
z}G_{k}[\varphi_{k}](e^{i\theta})\right|,\;\;\left|\frac{\partial}{\partial
\overline{z}}G_{k}[\varphi_{k}](e^{i\theta})\right|\Bigg\}\leq\nu_{k}^{\ast}(e^{i\theta}),$$
where $$\nu_{k}^{\ast}(e^{i\theta})=\begin{cases} \displaystyle
\frac{1}{4}\|\varphi_{1}\|_{\infty},
&\mbox{ if }\, k=1,\\
\displaystyle
\frac{1}{32}\left(\frac{3}{16}\right)^{k-2}\|\varphi_{k}\|_{\infty},
&\mbox{ if }\, 2\leq k\leq n-1.
\end{cases}$$
\end{lem}




\bpf
In order to prove the first statement of this Lemma, we
only need to prove the following inequality

$$\left|\frac{\partial}{\partial
z}G_{k}[\varphi_{k}](z)\right|\leq\nu_{k}(z)$$ because the proof of
the other one is similar, where $\nu_{k}$ is defined in the first
statement of this Lemma. For this, let \begin{eqnarray*} I_{2}(z)&=&
\int_{\mathbb{D}}\cdots\int_{\mathbb{D}}\left|\frac{\partial}{\partial
z}G(z,\xi_{1})\right|\cdots
|G(\xi_{k-1},\xi_{k})|\\
&&\times\left|\int_{0}^{2\pi}P(\xi_{k},e^{it})\varphi_{k}(e^{it})dt\right|d\sigma(\xi_{k})\cdots
d\sigma(\xi_{1}).
\end{eqnarray*}

\noindent $\mathbf{Case~ 1.}$  $k=1$.

Then, by \cite[Lemma 2.7]{K2}, we have
$$I_{2}(z)=
\int_{\mathbb{D}}\left|\frac{\partial}{\partial
z}G(z,\xi_{1})\right|\left|\int_{0}^{2\pi}P(\xi_{1},e^{it})\varphi_{1}(e^{it})dt\right|d\sigma(\xi_{1})
\leq\frac{1}{3}\|\varphi_{1}\|_{\infty},$$ which, together with
\cite[Proposition 2.4]{K2} (see also \cite{Tal}), gives that

\begin{eqnarray*}
\left|\frac{\partial}{\partial
z}G_{1}[\varphi_{1}](z)\right|&=&\left|\int_{\mathbb{D}}\frac{\partial}{\partial
z}G(z,\xi_{1})\left(\int_{0}^{2\pi}P(\xi_{1},e^{it})\varphi_{1}(e^{it})dt\right)d\sigma(\xi_{1})\right|\\
&\leq&I_{2}(z)\leq\frac{1}{3}\|\varphi_{1}\|_{\infty}.
\end{eqnarray*}

\noindent $\mathbf{Case~ 2.}$  $2\leq k=n-1$.

Let
$$\mathcal{A}_{1}=\int_{\mathbb{D}}\int_{\mathbb{D}}|G(\xi_{k-2},\xi_{k-1})|
|G(\xi_{k-1},\xi_{k})|\left|\int_{0}^{2\pi}P(\xi_{k},e^{it})\varphi_{k}(e^{it})dt\right|d\sigma(\xi_{k})d\sigma(\xi_{k-1}).$$

 It follows from Lemma \ref{lem-0.1} that

\begin{eqnarray*}
\mathcal{A}_{1}&\leq&\|\varphi_{k}\|_{\infty}\int_{\mathbb{D}}\int_{\mathbb{D}}|G(\xi_{k-2},\xi_{k-1})|
|G(\xi_{k-1},\xi_{k})|d\sigma(\xi_{k})d\sigma(\xi_{k-1})\\
&=&\frac{\|\varphi_{k}\|_{\infty}}{4}\int_{\mathbb{D}}|G(\xi_{k-2},\xi_{k-1})|(1-|\xi_{k-1}|^{2})d\sigma(\xi_{k-1})\\&\leq&
\frac{3(1-|\xi_{k-2}|^{2})}{64}\|\varphi_{k}\|_{\infty},
\end{eqnarray*}
which, together with Lemma \ref{lem-es-5}, implies that

\beq\label{eq-c11}
I_{2}(z)&\leq&\frac{3\|\varphi_{k}\|_{\infty}}{64}
\int_{\mathbb{D}}\cdots\int_{\mathbb{D}}\left|\frac{\partial}{\partial
z}G(z,\xi_{1})\right||G(\xi_{1},\xi_{2})|\cdots\\ \nonumber&&\times
|G(\xi_{k-3},\xi_{k-2})|(1-|\xi_{k-2}|^{2})d\sigma(\xi_{k-2})\cdots d\sigma(\xi_{1})\\
\nonumber&\leq&\frac{\|\varphi_{k}\|_{\infty}}{4}\left(\frac{3}{16}\right)^{k-2}\int_{\mathbb{D}}(1-|\xi_{1}|^{2})\left|\frac{\partial}{\partial
z}G(z,\xi_{1})\right| d\sigma(\xi_{1})\\ \nonumber
&=&\frac{\|\varphi_{k}\|_{\infty}}{16\pi}\left(\frac{3}{16}\right)^{k-2}\int_{\mathbb{D}}\frac{(1-|\xi_{1}|^{2})^{2}}{|1-z\overline{\xi}_{1}||z-\xi_{1}|}
d\sigma(\xi_{1})\\ \nonumber
&\leq&\|\varphi_{k}\|_{\infty}\left(\frac{3}{16}\right)^{k-2}\frac{(2-|z|^{2})}{30}.
\eeq





By (\ref{eq-c11}) and \cite[Proposition 2.4]{K2} (see also
\cite{Tal}), we conclude that

\begin{eqnarray*}
\left|\frac{\partial}{\partial
z}G_{k}[\varphi_{k}](z)\right|&=&\bigg|\int_{\mathbb{D}}\cdots\int_{\mathbb{D}}\frac{\partial}{\partial
z}G(z,\xi_{1})G(\xi_{1},\xi_{2})\cdots
G(\xi_{k-1},\xi_{k})\\
&&\times\left(\int_{0}^{2\pi}P(\xi_{k},e^{it})\varphi_{k}(e^{it})dt\right)d\sigma(\xi_{k})\cdots
d\sigma(\xi_{1})\bigg|\\
&\leq&I_{2}(z)\leq\|\varphi_{k}\|_{\infty}\left(\frac{3}{16}\right)^{k-2}\frac{(2-|z|^{2})}{30}.
\end{eqnarray*}

Now we prove the second statement of this Lemma. In order to show
this statement, we use the Vitali theorem (see \cite[Theorem
26.C]{Ha}) which asserts that if $\Omega$ is a measurable space with
finite measure $\mu$ and that
$\mathcal{F}_{n}:~\Omega\rightarrow\mathbb{C}$ is a sequence of
functions such that

$$\lim_{n\rightarrow\infty}\mathcal{F}_{n}(x)=\mathcal{F}(x)~\mbox{a.e.}~\mbox{and}~\sup_{n\geq1}\int_{\Omega}|\mathcal{F}_{n}|^{q}d\mu<\infty~\mbox{for some}~q>1,$$
then
$$\lim_{n\rightarrow\infty}\int_{\Omega}\mathcal{F}_{n}d\mu=\int_{\Omega}\mathcal{F}d\mu.$$

Let

\begin{eqnarray*}
\mathcal{A}_{2}&=&\int_{\mathbb{D}}\bigg(\bigg|\frac{\partial}{\partial
z}G(z,\xi_{1})\bigg|^{\frac{3}{2}}\bigg|\int_{\mathbb{D}}\bigg(G(\xi_{1},\xi_{2})\cdots\int_{\mathbb{D}}
\bigg(G(\xi_{k-1},\xi_{k})\\
&&\times\bigg(\int_{0}^{2\pi}P(\xi_{k},e^{it})\varphi_{k}(e^{it})dt\bigg)\bigg)d\sigma(\xi_{k})\cdots\bigg)d\sigma(\xi_{2})
\bigg|^{\frac{3}{2}}\bigg)d\sigma(\xi_{1}),
\end{eqnarray*} where $k\in\{1,\ldots,n-1\}$.
In order to estimate $\mathcal{A}_{2}$, we let
\be\label{eq-kk-6}\eta_{1}=\frac{z-\xi_{1}}{1-\overline{z}\xi_{1}}=r_{1}
e^{i\theta_{1}},\ee where $r_{1}\in[0,1)$ and $t\in[0,2\pi]$. Since
\begin{eqnarray*}
\int_{0}^{1}r_{1}^{2j-\frac{1}{2}}(1-r_{1}^{2})^{\frac{3}{2}}dr_{1}&=&\frac{3}{\big(2j+\frac{5}{2}\big)\big(2j+\frac{1}{2}\big)}\int_{0}^{1}\frac{x^{2j+\frac{5}{2}}}{(1-x^{2})^{\frac{1}{2}}}dx\\
&=&\frac{3}{\big(2j+\frac{5}{2}\big)\big(2j+\frac{1}{2}\big)}\int_{0}^{\frac{\pi}{2}}\left(\sin
t\right)^{2j+\frac{5}{2}}dt\\
&\leq&\frac{3}{\big(2j+\frac{5}{2}\big)\big(2j+\frac{1}{2}\big)}\int_{0}^{\frac{\pi}{2}}\left(\sin
t\right)^{2j+2}dt\\
&=&\frac{3\pi}{2\big(2j+\frac{5}{2}\big)\big(2j+\frac{1}{2}\big)}\cdot\frac{(2j+1)!!}{(2j+2)!!},
\end{eqnarray*}
we see that

\begin{eqnarray*}
\sum_{j=0}^{\infty}(j+1)^{2}|z|^{2j}\int_{0}^{1}r_{1}^{2j-\frac{1}{2}}(1-r_{1}^{2})^{\frac{3}{2}}dr_{1}&\leq&\frac{6\pi}{5}\sum_{j=0}^{\infty}\frac{(2j+1)!!}{(2j+2)!!}|z|^{2j}\\
&<&\frac{6\pi}{5}\left(1+\sum_{j=1}^{\infty}\frac{(2j-1)!!}{(2j)!!}|z|^{2j}\right)\\
&=&\frac{6\pi}{5}\frac{1}{(1-|z|^{2})^{\frac{1}{2}}},
\end{eqnarray*}
which, together with  (\ref{eq-kk-6}) and Theorem \Ref{Thm-1}, yield
that

\beq\label{kk-7}\int_{\mathbb{D}}\bigg|\frac{\partial}{\partial
z}G(z,\xi_{1})\bigg|^{\frac{3}{2}}d\sigma(\xi_{1})&=&\frac{1}{2\pi}\int_{\mathbb{D}}
\frac{(1-|\eta_{1}|^{2})^{\frac{3}{2}}(1-|z|^{2})^{\frac{1}{2}}}{|\eta_{1}|^{\frac{3}{2}}|1-\overline{z}\eta_{1}|^{4}}d\sigma(\eta_{1})\\
\nonumber
&=&(1-|z|^{2})^{\frac{1}{2}}\sum_{j=0}^{\infty}(j+1)^{2}|z|^{2j}\\
\nonumber
&&\times\int_{0}^{1}r_{1}^{2j-\frac{1}{2}}(1-r_{1}^{2})^{\frac{3}{2}}dr_{1}
<\frac{6\pi}{5}. \eeq
It follows from   (\ref{kk-1}) and (\ref{kk-7}) that

\begin{eqnarray*}
\mathcal{A}_{2}&\leq&\|\varphi_{k}\|_{\infty}^{\frac{3}{2}}\int_{\mathbb{D}}\bigg(\bigg|\frac{\partial}{\partial
z}G(z,\xi_{1})\bigg|^{\frac{3}{2}}\\
&&\times\bigg(\int_{\mathbb{D}}\bigg(|G(\xi_{1},\xi_{2})|\cdots\int_{\mathbb{D}}
|G(\xi_{k-1},\xi_{k})|d\sigma(\xi_{k})\cdots\bigg)d\sigma(\xi_{2})
\bigg)^{\frac{3}{2}}\bigg)d\sigma(\xi_{1})\\
&\leq&\|\varphi_{k}\|_{\infty}^{\frac{3}{2}}\left(\frac{1}{4}\right)^{\frac{3(k-1)}{2}}\int_{\mathbb{D}}\bigg|\frac{\partial}{\partial
z}G(z,\xi_{1})\bigg|^{\frac{3}{2}}d\sigma(\xi_{1})\\
&\leq&\frac{5\|\varphi_{k}\|_{\infty}^{\frac{3}{2}}\pi}{6}\left(\frac{1}{8}\right)^{k-1}<\infty.
\end{eqnarray*}
Therefore, by the Vitali theorem,  we conclude that
$\frac{\partial}{\partial z}G_{k}[\varphi_{k}]$ has continuous
extension to the boundary, and further, by \cite[Lemma 2.7]{K2},

$$\left|\frac{\partial}{\partial
z}G_{1}[\varphi_{1}](e^{i\theta})\right|\leq\frac{1}{4}\|\varphi_{1}\|_{\infty},$$
where $\theta\in[0,2\pi]$.

For $2\leq k=n-1$, by Lemmas \ref{lem-0.1} and \ref{lem-y1},  we
have

\begin{eqnarray*} \left|\frac{\partial}{\partial
z}G_{k}[\varphi_{k}](e^{i\theta})\right|&=&\bigg|\int_{\mathbb{D}}\cdots\int_{\mathbb{D}}\frac{\partial}{\partial
z}G(e^{i\theta},\xi_{1})G(\xi_{1},\xi_{2})\cdots
G(\xi_{k-1},\xi_{k})\\
&&\times\left(\int_{0}^{2\pi}P(\xi_{k},e^{it})\varphi_{k}(e^{it})dt\right)d\sigma(\xi_{k})\cdots
d\sigma(\xi_{1})\bigg|\\
&\leq&\frac{1}{4}\left(\frac{3}{16}\right)^{k-2}\|\varphi_{k}\|_{\infty}\int_{\mathbb{D}}\left|\frac{\partial}{\partial
z}G(e^{i\theta},\xi_{1})\right|(1-|\xi_{1}|^{2})d\sigma(\xi_{1})\\
&=&\frac{1}{32}\left(\frac{3}{16}\right)^{k-2}\|\varphi_{k}\|_{\infty}.
\end{eqnarray*}

Similarly, we can show that $\frac{\partial}{\partial
\overline{z}}G_{k}[\varphi_{k}]$ has continuous extension to the
boundary, and  for $\theta\in[0,2\pi]$,
$$\left|\frac{\partial}{\partial
\overline{z}}G_{k}[\varphi_{k}](e^{i\theta})\right|\leq\nu_{k}^{\ast}(e^{i\theta}),$$
where $k\in\{1,\ldots,n-1\}$.
 The proof of this lemma is complete.
 \epf \medskip


\begin{lem}\label{lem-2} Suppose $\varphi_{n}\in\mathcal{C}(\overline{\mathbb{D}})$ and $G_{n}[\varphi_{n}]$ is defined in \eqref{eq-ch-3.2}.
Then, the following statements hold:

\noindent $(1)$
 For
$z\in\mathbb{D}$, $$\max\Bigg\{\left|\frac{\partial}{\partial z}
G_{n}[\varphi_{n}](z)\right|,~\left|\frac{\partial}{\partial
\overline{z}}
G_{n}[\varphi_{n}](z)\right|\Bigg\}\leq\|\varphi_{n}\|_{\infty}\left(\frac{3}{16}\right)^{n-2}\frac{(2-|z|^{2})}{30}.$$

\noindent $(2)$  Both $\frac{\partial}{\partial z}
G_{n}[\varphi_{n}]$ and $ \frac{\partial}{\partial \overline{z}}
G_{n}[\varphi_{n}]$ have continuous extensions to the boundary, and
further, for $\theta\in[0,2\pi]$,

$$\max\Bigg\{\left|\frac{\partial}{\partial z}
G_{n}[\varphi_{n}](e^{i\theta})\right|,~\left|\frac{\partial}{\partial
\overline{z}}
G_{n}[\varphi_{n}](e^{i\theta})\right|\Bigg\}\leq\frac{1}{32}\left(\frac{3}{16}\right)^{n-2}\|\varphi_{n}\|_{\infty}.$$

\end{lem}

\bpf To prove the first statement, we only need to prove the
inequality:

$$\left|\frac{\partial}{\partial z}
G_{n}[\varphi_{n}](z)\right|\leq\|\varphi_{n}\|_{\infty}\left(\frac{3}{16}\right)^{n-2}\frac{(2-|z|^{2})}{30},~z\in\mathbb{D},$$
because the proof to the other one is similar. For this,  let

\begin{eqnarray*}
 I_{3}(z)&=&
\int_{\mathbb{D}}\cdots\int_{\mathbb{D}}\left|\frac{\partial}{\partial z}G(z,\zeta_{1})\right||G(\zeta_{1},\zeta_{2})|\cdots |G(\zeta_{n-2},\zeta_{n-1})|\\
\nonumber
&&\times\left|\int_{\mathbb{D}}G(\zeta_{n-1},\zeta_{n})\varphi_{n}(\zeta_{n})d\sigma(\zeta_{n})\right|d\sigma(\zeta_{n-1})\cdots
d\sigma(\zeta_{1}).
\end{eqnarray*}

By calculation, we have
\beq\label{eq-ch-0.1}\left|\int_{\mathbb{D}}G(\zeta_{n-1},\zeta_{n})\varphi_{n}(\zeta_{n})d\sigma(\zeta_{n})\right|
&\leq&\|\varphi_{n}\|_{\infty}\int_{\mathbb{D}}|G(\zeta_{n-1},\zeta_{n})|d\sigma(\zeta_{n})\\
\nonumber
&=&\frac{\|\varphi_{n}\|_{\infty}(1-|\zeta_{n-1}|^{2})}{4}.\eeq

By (\ref{eq-ch-0.1}), Lemmas \ref{lem-0.1} and \ref{lem-es-5}, we
see that

\begin{eqnarray*}
 I_{3}(z)&\leq&\frac{\|\varphi_{n}\|_{\infty}}{4}
\int_{\mathbb{D}}\cdots\int_{\mathbb{D}}\left|\frac{\partial}{\partial z}G(z,\zeta_{1})\right||G(\zeta_{1},\zeta_{2})|\cdots |G(\zeta_{n-2},\zeta_{n-1})|\\
\nonumber &&\times(1-|\zeta_{n-1}|^{2})d\sigma(\zeta_{n-1})\cdots
d\sigma(\zeta_{1})\\
&\leq&\frac{\|\varphi_{n}\|_{\infty}}{16\pi}\left(\frac{3}{16}\right)^{n-2}\int_{\mathbb{D}}\frac{(1-|\zeta_{1}|^{2})^{2}}{|1-z\overline{\zeta}_{1}||z-\zeta_{1}|}
d\sigma(\zeta_{1})\\
&\leq&\|\varphi_{n}\|_{\infty}\left(\frac{3}{16}\right)^{n-2}\frac{(2-|z|^{2})}{30},
\end{eqnarray*}
which, together with \cite[Proposition 2.4]{K2} (see also
\cite{Tal}), yields that \beq\label{eq-ch-0.2}
\left|\frac{\partial}{\partial z}G_{n}[\varphi_{n}](z)\right|&=&
\bigg|\int_{\mathbb{D}}\cdots\int_{\mathbb{D}}\frac{\partial}{\partial
z}G(z,\zeta_{1})\cdots G(\zeta_{n-2},\zeta_{n-1})\\
\nonumber
&&\times\left[\int_{\mathbb{D}}G(\zeta_{n-1},\zeta_{n})\varphi_{n}(\zeta_{n})d\sigma(\zeta_{n})\right]d\sigma(\zeta_{n-1})\cdots
d\sigma(\zeta_{1})\bigg|\\
\nonumber&\leq&I_{3}\leq\|\varphi_{n}\|_{\infty}\left(\frac{3}{16}\right)^{n-2}\frac{(2-|z|^{2})}{30}.\eeq

Next, we prove the second part of this Lemma. Set

\begin{eqnarray*}
\mathcal{A}_{3}&=&\int_{\mathbb{D}}\bigg(\bigg|\frac{\partial}{\partial
z}G(z,\zeta_{1})\bigg|^{\frac{3}{2}}\bigg|\int_{\mathbb{D}}\bigg(G(\zeta_{1},\zeta_{2})\cdots\int_{\mathbb{D}}
\bigg(G(\zeta_{n-2},\zeta_{n-1})\\
&&\times\bigg(\int_{\mathbb{D}}G(\zeta_{n-1},\zeta_{n})\varphi_{n}(\zeta_{n})d\sigma(\zeta_{n})\bigg)\bigg)d\sigma(\zeta_{n-1})\cdots\bigg)d\sigma(\zeta_{2})
\bigg|^{\frac{3}{2}}\bigg)d\sigma(\zeta_{1}).
\end{eqnarray*}
Then by (\ref{kk-7}) and Lemma \ref{lem-0.1} (\ref{kk-1}), we get
\begin{eqnarray*}
 \mathcal{A}_{3}&\leq&
\frac{\|\varphi_{n}\|_{\infty}^{\frac{3}{2}}}{8}\int_{\mathbb{D}}\bigg(\bigg|\frac{\partial}{\partial
z}G(z,\zeta_{1})\bigg|^{\frac{3}{2}}\bigg|\int_{\mathbb{D}}\bigg(G(\zeta_{1},\zeta_{2})\cdots\int_{\mathbb{D}}
\bigg(G(\zeta_{n-2},\zeta_{n-1})\\
&&\times\big(1-|\zeta_{n-1}|^{2}\big)\bigg)d\sigma(\zeta_{n-1})\cdots\bigg)d\sigma(\zeta_{2})
\bigg|^{\frac{3}{2}}\bigg)d\sigma(\zeta_{1})\\
&\leq&\|\varphi_{n}\|_{\infty}^{\frac{3}{2}}\left(\frac{1}{8}\right)^{n-1}\int_{\mathbb{D}}\bigg|\frac{\partial}{\partial
z}G(z,\zeta_{1})\bigg|^{\frac{3}{2}}d\sigma(\zeta_{1})\\
&\leq&\frac{6\pi}{5}\left(\frac{1}{8}\right)^{n-1}\|\varphi_{n}\|_{\infty}^{\frac{3}{2}}<\infty.
\end{eqnarray*}
Hence, by the Vitali theorem,  we see that $\frac{\partial}{\partial
z}G_{n}[\varphi_{n}]$ has continuous extension to the boundary, and
further, by Lemmas \ref{lem-0.1} and \ref{lem-y1},  we have

\begin{eqnarray*} \left|\frac{\partial}{\partial
z}G_{n}[\varphi_{n}](e^{i\theta})\right|&=&\bigg|\int_{\mathbb{D}}\cdots\int_{\mathbb{D}}\frac{\partial}{\partial
z}G(e^{i\theta},\zeta_{1})G(\zeta_{1},\zeta_{2})\cdots
G(\zeta_{n-1},\zeta_{n})\\
&&\times\left(\int_{\mathbb{D}}G(\zeta_{n-1},\zeta_{n})\varphi_{n}(\zeta_{n})d\sigma(\zeta_{n})\right)d\sigma(\zeta_{n-1})\cdots
d\sigma(\zeta_{1})\bigg|\\
&\leq&\frac{1}{4}\left(\frac{3}{16}\right)^{n-2}\|\varphi_{n}\|_{\infty}\int_{\mathbb{D}}\left|\frac{\partial}{\partial
z}G(e^{i\theta},\zeta_{1})\right|(1-|\zeta_{1}|^{2})d\sigma(\zeta_{1})\\
&=&\frac{1}{32}\left(\frac{3}{16}\right)^{n-2}\|\varphi_{n}\|_{\infty},
\end{eqnarray*}
where $\theta\in[0,2\pi]$.

Similarly, we can prove that $\frac{\partial}{\partial
\overline{z}}G_{n}[\varphi_{n}]$ has continuous extension to the
boundary, and  for $\theta\in[0,2\pi]$,
$$\left|\frac{\partial}{\partial
\overline{z}}G_{n}[\varphi_{n}](e^{i\theta})\right|\leq\frac{\|\varphi_{n}\|_{\infty}}{32}\left(\frac{3}{16}\right)^{n-2}.$$
The proof of this lemma is finished.
 \epf

\begin{lem}\label{lem-main}
For  $\varphi_{k}\in\mathcal{C}(\mathbb{T})$ and
$\varphi_{n}\in\mathcal{C}(\overline{\mathbb{D}})$, suppose that $f$
is a sense-preserving
 homeomorphism from $\overline{\mathbb{D}}$ onto  itself
satisfying {\rm (\ref{eq-ch-1})} and the boundary conditions
$\Delta^{n-1}f=\varphi_{n-1},~\ldots,~\Delta^{1}f=\varphi_{1}$ on
$\mathbb{T}$, and suppose that $f$ is Lipschitz continuous in
$\mathbb{D}$, where $n\geq2$ and $k\in\{1,\ldots,n-1\}$. Then, for
almost every $e^{i\theta}\in\mathbb{T}$, the following limits exist:
\be \label{qw-1} D_{f}(e^{i\theta}):=\lim_{z\rightarrow
e^{i\theta},z\in\mathbb{D}}D_{f}(z)\;\;\mbox{and}\;\;
J_{f}(e^{i\theta}):=\lim_{z\rightarrow
e^{i\theta},z\in\mathbb{D}}J_{f}(z).\ee Further, we have
\beq\label{eq-sh-1}J_{f}(e^{i\theta})&\leq&\frac{\gamma'(\theta)}{2\pi}\int_{0}^{2\pi}\frac{|f(e^{it})-f(e^{i\theta})|^{2}}{|e^{it}-e^{i\theta}|^{2}}dt
+\frac{\gamma'(\theta)\|\varphi_{1}\|_{\infty}}{2}\\
\nonumber&&+\gamma'(\theta)\sum_{k=2}^{n}\frac{\|\varphi_{k}\|_{\infty}}{16}\left(\frac{3}{16}\right)^{k-2}
\eeq and
\beq\label{eq-sh-1.1}J_{f}(e^{i\theta})&\geq&\frac{\gamma'(\theta)}{2\pi}\int_{0}^{2\pi}\frac{|f(e^{it})-f(e^{i\theta})|^{2}}{|e^{it}-e^{i\theta}|^{2}}dt
-\frac{\gamma'(\theta)\|\varphi_{1}\|_{\infty}}{2}\\
\nonumber&&-\gamma'(\theta)\sum_{k=2}^{n}\frac{\|\varphi_{k}\|_{\infty}}{16}\left(\frac{3}{16}\right)^{k-2},
\eeq where $f(e^{i\theta})=e^{i\gamma(\theta)}$ and $\gamma(\theta)$
is a real-valued function in $[0,2\pi]$.
\end{lem}





\subsection*{Proof of Lemma \ref{lem-main}} We first prove the existence of the two limits in (\ref{qw-1}). By Lemmas \ref{lem-1} and
\ref{lem-2}, we get that for any $e^{i\theta}\in\mathbb{D}$,

\be\label{eq-kk-1}\lim_{z\rightarrow
e^{i\theta},z\in\mathbb{D}}D_{G_{k}[\varphi_{k}]}(z)=D_{G_{k}[\varphi_{k}]}(e^{i\theta}),\ee
where $k\in\{1,2,\ldots,n\}$.

Again, by Lemmas  \ref{lem-1} and \ref{lem-2}, we know that
$\|D_{G_{k}[\varphi_{k}]}\|$ is bounded, which implies the Lipschitz
continuity of $G_{k}[\varphi_{k}]$  in $\mathbb{D}$. Since $f$ is
Lipschitz continuous in $\mathbb{D}$, we see that $\|D_{f}\|$ is
bounded in $\mathbb{D}$. Thus, it follows from (\ref{eq-ch-3.0})
that $$P[\varphi_{0}]=f-\sum_{k=1}^{n}(-1)^{k}G_{k}[\varphi_{k}]$$
is also Lipschitz continuous in $\mathbb{D}$, where
$\varphi_{0}=f|_{\mathbb{T}}$. Now, we conclude from \cite[Lemma
2.1]{K2} that
 for almost every $e^{i\theta}\in\mathbb{T}$,
$$\lim_{z\rightarrow
e^{i\theta},z\in\mathbb{D}}D_{P[\varphi_{0}]}(z)$$ does exist,
which, together with (\ref{eq-ch-3.0}) and (\ref{eq-kk-1}),
guarantees that for almost every $\theta\in[0,2\pi],$
$$\lim_{z\rightarrow
e^{i\theta},z\in\mathbb{D}}D_{f}(z)$$ also exists.

Since $$J_{f}(z)=\det
D_{f}(z),$$ obviously, we see that  $$\lim_{z\rightarrow
e^{i\theta},z\in\mathbb{D}}J_{f}(z)$$ exists for almost every
$\theta\in[0,2\pi].$

Next, we demonstrate the estimates in (\ref{eq-sh-1}) and  (\ref{eq-sh-1.1}). For convenience, in the rest of the proof of the lemma, let
$$D_{f}(e^{i\theta})=\lim_{z\rightarrow e^{i\theta},z\in\mathbb{D}}D_{f}(z)\;\;\mbox{and}\;\; J_{f}(e^{i\theta})=\lim_{z\rightarrow
e^{i\theta},z\in\mathbb{D}}J_{f}(z).$$


By Lebesgue Dominated Convergence Theorem, the boundedness of
$\|D_{f}\|$, and by letting $z=re^{i\theta}\in\mathbb{D}$, we see
that for any fixed $\theta_{0}\in[0,2\pi]$, \beq\label{eq-d1}
f(e^{i\theta})&=&\lim_{r\rightarrow1^{-}}f(re^{i\theta})=\lim_{r\rightarrow1^{-}}\int_{\theta_{0}}^{\theta}\frac{\partial}{\partial
t}f(re^{it})dt+f(e^{i\theta_{0}})\\ \nonumber
&=&\int_{\theta_{0}}^{\theta}\lim_{r\rightarrow1^{-}}\left[ir\big(f_{z}(re^{it})e^{it}-f_{\overline{z}}(re^{it})e^{-it}\big)\right]dt+f(e^{i\theta_{0}}),
\eeq which implies that $f(e^{i\theta})$ is absolutely continuous.
Let $\gamma(\theta)$ be a real-valued function in $[0,2\pi]$ such
that
$$e^{i\gamma(\theta)}=f(e^{i\theta}).$$
Then,
\be\label{eq-r1}f'(e^{i\theta})=i\gamma'(\theta)e^{i\gamma(\theta)}\ee
holds almost everywhere in $[0,2\pi].$

Since
$$J_{f}(re^{i\theta})=|f_{z}(re^{i\theta})|^{2}-|f_{\overline{z}}(re^{i\theta})|^{2}=
-\mbox{Re}\left(\overline{\frac{\partial f}{\partial
r}}\frac{i}{r}\frac{\partial f}{\partial \theta}\right),$$ we infer
from (\ref{eq-r1}) that

\beq\label{sh-2}
J_{f}(e^{i\theta})&=&\lim_{r\rightarrow1^{-}}J_{f}(re^{i\theta})=-\lim_{r\rightarrow1^{-}}\mbox{Re}\left(\overline{\frac{\partial
f}{\partial r}}\frac{i}{r}\frac{\partial f}{\partial
\theta}\right)=\chi_{0}-\sum_{k=1}^{n}(-1)^{k}\chi_{k},\eeq where
$$\chi_{0}=\lim_{r\rightarrow1^{-}}\mbox{Re}\left(\frac{\overline{f(e^{i\theta})-P[\varphi_{0}](re^{i\theta})}}{1-r}\cdot\gamma'(\theta)f(e^{i\theta})\right)$$
and
$$\chi_{k}=\lim_{r\rightarrow1^{-}}\mbox{Re}\left(\frac{\overline{G_{k}[\varphi_{k}](re^{i\theta})}}{1-r}\cdot\gamma'(\theta)f(e^{i\theta})\right).$$

Now, we are going to prove (\ref{eq-sh-1}) and (\ref{eq-sh-1.1}) by
estimating the quantities $|\chi_{0}|$ and $|\chi_{k}|,$
respectively, where $k\in\{1,2,\ldots,n\}.$ We start with the
estimate of $|\chi_{0}|$. Since
$$
\mbox{Re}\langle
f(e^{i\theta}),f(e^{i\theta})-f(e^{it})\rangle=\mbox{Re}\big[f(e^{i\theta})(\overline{f(e^{i\theta})-f(e^{it})})\big]
=\frac{1}{2}|f(e^{it})-f(e^{i\theta})|^{2}
 $$ and
$$\chi_{0}=\lim_{r\rightarrow1^{-}}\mbox{Re}\left(\frac{1}{2\pi}\int_{0}^{2\pi}\frac{1+r}{|1-re^{i(\theta-t)}|^{2}}\langle
\gamma'(\theta)f(e^{i\theta}),f(e^{i\theta})-f(e^{it})\rangle\,dt\right),$$
where $\langle\cdot,\cdot\rangle$ denotes the inner product, it
follows that \be\label{eq-s3}
\chi_{0}=\gamma'(\theta)\frac{1}{2\pi}\int_{0}^{2\pi}\frac{|f(e^{it})-f(e^{i\theta})|^{2}}{|e^{it}-e^{i\theta}|^{2}}dt.
\ee

Next, we estimate $|\chi_{k}|$ for $k\in\{1,2,\ldots,n\}.$ Since
\be\label{eqr1}\lim_{z\rightarrow
e^{i\theta},z\in\mathbb{D}}\frac{G(z,\xi_{1})}{1-|z|}=\lim_{z\rightarrow
e^{i\theta},z\in\mathbb{D}}\frac{G(z,\xi_{1})-G(e^{i\theta},\xi_{1})}{1-|z|}=P(\xi_{1},e^{i\theta}),\ee
 we deduce that

\beq\label{th-1}\lim_{z\rightarrow
e^{i\theta},z\in\mathbb{D}}\frac{G_{k}[\varphi_{k}](z)}{1-|z|}&=&\lim_{z\rightarrow
e^{i\theta},z\in\mathbb{D}}\int_{\mathbb{D}}\cdots\int_{\mathbb{D}}\frac{G(z,\xi_{1})}{1-|z|}G(\xi_{1},\xi_{2})\cdots
G(\xi_{k-1},\xi_{k})\\ \nonumber
&&\times\left[\int_{0}^{2\pi}P(\xi_{k},e^{it})\varphi_{k}(e^{it})dt\right]d\sigma(\xi_{k})\cdots
d\sigma(\xi_{1})\\ \nonumber &=&
\int_{\mathbb{D}}\cdots\int_{\mathbb{D}}P(\xi_{1},e^{i\theta})G(\xi_{1},\xi_{2})\cdots
G(\xi_{k-1},\xi_{k})\\ \nonumber
&&\times\left[\int_{0}^{2\pi}P(\xi_{k},e^{it})\varphi_{k}(e^{it})dt\right]d\sigma(\xi_{k})\cdots
d\sigma(\xi_{1}),\eeq where $k\in\{1,2,\ldots,n-1\}.$

\noindent $\mathbf{Case~ 1.}$  $2\leq k\leq n-1.$

It follows from (\ref{th-1}), Lemmas \ref{lem-0.1} and \ref{lem-y1}
that

\beq\label{eqr2}
 |\chi_{k}|&\leq&
\gamma'(\theta)\|\varphi_{k}\|_{\infty}\int_{\mathbb{D}}\cdots\int_{\mathbb{D}}P(\xi_{1},e^{i\theta})|G(\xi_{1},\xi_{2})|\cdots\\
\nonumber &&|G(\xi_{k-1},\xi_{k})|d\sigma(\xi_{k})\cdots
d\sigma(\xi_{1})\\ \nonumber
&\leq&\frac{\gamma'(\theta)\|\varphi_{k}\|_{\infty}}{4}\left(\frac{3}{16}\right)^{k-2}\int_{\mathbb{D}}P(\xi_{1},e^{i\theta})(1-|\xi_{1}|^{2})d\sigma(\xi_{1})\\
\nonumber
&=&\frac{\gamma'(\theta)\|\varphi_{k}\|_{\infty}}{16}\left(\frac{3}{16}\right)^{k-2}.
\eeq

\noindent $\mathbf{Case~ 2.}$  $k=1.$

By  (\ref{th-1}), we get

\be\label{eq-r2.0}
 |\chi_{1}|\leq\gamma'(\theta)\|\varphi_{1}\|_{\infty}\int_{\mathbb{D}}P(\xi_{1},e^{i\theta})d\sigma(\xi_{1})=\frac{\gamma'(\theta)\|\varphi_{1}\|_{\infty}}{2}.
\ee

At last, we estimate $ |\chi_{n}|$. By (\ref{eqr1}), Lemmas
\ref{lem-0.1} and \ref{lem-y1}, we obtain

\beq\label{eqr3}
 |\chi_{n}|&\leq&\gamma'(\theta)\int_{\mathbb{D}}\cdots\int_{\mathbb{D}}P(\zeta_{1},e^{i\theta})|G(\zeta_{1},\zeta_{2})|\cdots
|G(\zeta_{n-2},\zeta_{n-1})|\\ \nonumber
&&\times\left|\int_{\mathbb{D}}G(\zeta_{n-1},\zeta_{n})
\varphi_{n}(\zeta_{n})d\sigma(\zeta_{n})\right|d\sigma(\zeta_{n-1})\cdots
d\sigma(\zeta_{1})\\ \nonumber &\leq&
\gamma'(\theta)\|\varphi_{n}\|_{\infty}\int_{\mathbb{D}}\cdots\int_{\mathbb{D}}P(\zeta_{1},e^{i\theta})|G(\zeta_{1},\zeta_{2})|\cdots\\
\nonumber &&|G(\zeta_{n-1},\zeta_{n})|d\sigma(\zeta_{n})\cdots
d\sigma(\zeta_{1})\\ \nonumber
&\leq&\frac{\gamma'(\theta)\|\varphi_{n}\|_{\infty}}{4}\left(\frac{3}{16}\right)^{n-2}\int_{\mathbb{D}}P(\zeta_{1},e^{i\theta})(1-|\zeta_{1}|^{2})d\sigma(\zeta_{1})\\
\nonumber
&=&\frac{\gamma'(\theta)\|\varphi_{n}\|_{\infty}}{16}\left(\frac{3}{16}\right)^{n-2}.
\eeq

Hence, (\ref{eq-sh-1}) and (\ref{eq-sh-1.1}) follow from the
inequalities (\ref{eq-s3}), (\ref{eqr2}) and (\ref{eqr3}) along with
the following chain of inequalities:
$$|\chi_{0}|-\sum_{k=1}^{n}|\chi_{k}|\leq J_{f}(e^{i\theta})\leq \sum_{k=0}^{n}|\chi_{k}|.$$
The proof of the lemma is complete. \qed
\medskip

\begin{Thm}{\rm (\cite[Theorem 3.4]{K-M})}\label{KM-P}
Suppose that $f$ is a quasiconformal $\mathcal{C}^{2}$ diffeomorphism from
the plane domain $\Omega$ with $\mathcal{C}^{1,\alpha}$ compact
boundary onto the  plane domain $\Omega^{\ast}$ with
$\mathcal{C}^{2,\alpha}$ compact boundary. If there exist constants
$a_{1}$ and $b_{1}$ such that
$$|\Delta f(z)|\leq a_{1}\|D_{f}(z)\|^{2}+b_{1}$$ in $\Omega$, then $f$ has
bounded partial derivatives. In particular, it is a
Lipschitz mapping in $\Omega$.
\end{Thm}

\begin{Thm}{\rm (\cite[Theorem 2.2]{PS-2005})}\label{PS}
Given $K\geq1$, let $f$ be a $K$-quasiconformal and harmonic
self-mapping of $\mathbb{D}$ satisfying $f(0)=0$. Then, for
$z\in\mathbb{D}$,

$$|f_{z}(z)|\geq\frac{K+1}{2K}\max\left\{\frac{2}{\pi},L_{K}\right\},$$
where $L_{K}$ is defined in \cite[Lemma 1.4]{PS-2005} and
$\lim_{K\rightarrow1}L_{K}=1.$
\end{Thm}

\section{The proof of Theorem \ref{thm-1.12}}\label{sec-3}
We first prove part $(a)$.

\bst\label{bst-4.1} The co-Lipschitz continuity of $f$.\est

Now we begin to prove the co-Lipschitz continuity of $f$. Since $f$
is a $(K,K')$-quasiconformal mapping, by \cite[Lemma 4.2]{K8}, we
see that, for $z\in\mathbb{D}$,
$$\|D_{f}(z)\|\leq K\lambda(D_{f}(z))+\sqrt{K'},$$
which implies that

\be\label{KM-1}|f_{\overline{z}}(z)|\leq\frac{(K-1)}{K+1}|f_{z}(z)|+\frac{\sqrt{K'}}{K+1}.\ee
By (\ref{eq-ch-3.0}), we have
$$|f_{\overline{z}}(z)|=\left|P[\varphi_{0}]_{\overline{z}}(z)+\sum_{k=1}^{n}(-1)^{k}G_{k}[\varphi_{k}]_{\overline{z}}(z)\right|\geq
\left|P[\varphi_{0}]_{\overline{z}}(z)\right|-\sum_{k=1}^{n}\left|G_{k}[\varphi_{k}]_{\overline{z}}(z)\right|$$
and
$$|f_{z}(z)|=\left|P[\varphi_{0}]_{z}(z)+\sum_{k=1}^{n}(-1)^{k}G_{k}[\varphi_{k}]_{z}(z)\right|\leq
\left|P[\varphi_{0}]_{z}(z)\right|+\sum_{k=1}^{n}\left|G_{k}[\varphi_{k}]_{z}(z)\right|,$$
which, together with (\ref{KM-1}), yield that

\be\label{KM-2}|P[\varphi_{0}]_{\overline{z}}(z)|\leq\frac{(K-1)}{K+1}|P[\varphi_{0}]_{z}(z)|+\Lambda(z),\ee
where
$$\Lambda(z)=\frac{(K-1)}{K+1}\sum_{k=1}^{n}\left|G_{k}[\varphi_{k}]_{z}(z)\right|+\sum_{k=1}^{n}\left|G_{k}[\varphi_{k}]_{\overline{z}}(z)\right|+\frac{\sqrt{K'}}{1+K}.$$
By Lemmas \ref{lem-1} and \ref{lem-2}, we have

\be\label{KM-3}\Lambda(z)\leq\frac{2K}{K+1}\mathcal{H}(\varphi_{1},\cdots,\varphi_{n})+\frac{\sqrt{K'}}{1+K},
\ee where
$$\mathcal{H}(\varphi_{1},\cdots,\varphi_{n})=\left(\frac{1}{3}\|\varphi_{1}\|_{\infty}+\frac{1}{15}\sum_{k=2}^{n}\left(\frac{3}{16}\right)^{k-2}\|\varphi_{k}\|_{\infty}\right).$$

Since $\varphi_{0}$ is a sense-preserving homeomorphic self-mapping
of $\mathbb{T}$, by the Choquet-Rad\'o-Kneser theorem (see
\cite{Cho}), we see that $P[\varphi_{0}]$ is a harmonic
diffeomorphism of $\mathbb{D}$ onto itself. Then, by \cite[Lemma
2.1]{ZT}, we obtain
\be\label{KM-4}\frac{|P[\varphi_{0}]_{z}(z)|}{\frac{1}{\pi}-\frac{|P[\varphi_{0}](0)|}{2}}\geq1.\ee
It follows from (\ref{KM-2}), (\ref{KM-3}) and (\ref{KM-4}) that

$$|P[\varphi_{0}]_{\overline{z}}(z)|\leq\frac{(K-1)}{K+1}|P[\varphi_{0}]_{z}(z)|+
\left(\frac{2K}{K+1}\mathcal{H}(\varphi_{1},\cdots,\varphi_{n})+\frac{\sqrt{K'}}{K+1}\right)\frac{|P[\varphi_{0}]_{z}(z)|}{\left(\frac{1}{\pi}-\frac{|P[\varphi_{0}](0)|}{2}\right)},$$
which, together with the assumptions, gives that

\be\label{KM-5}\frac{|P[\varphi_{0}]_{\overline{z}}(z)|}{|P[\varphi_{0}]_{z}(z)|}\leq
\frac{K-1}{K+1}+\frac{2K\mathcal{H}(\varphi_{1},\cdots,\varphi_{n})+\sqrt{K'}}{(K+1)\left(\frac{1}{\pi}-\frac{|P[\varphi_{0}](0)|}{2}\right)}<1.\ee
Then $P[\varphi_{0}]$ is a $K^{\ast}$-quasiconformal mapping in
$\mathbb{D}$, where

\be\label{KM-5.1}K^{\ast}=\frac{K\left(\frac{2}{\pi}-|P[\varphi_{0}](0)|\right)+2K\mathcal{H}(\varphi_{1},\cdots,\varphi_{n})+
\sqrt{K'}}{\frac{2}{\pi}-|P[\varphi_{0}](0)|-2K\mathcal{H}(\varphi_{1},\cdots,\varphi_{n})-\sqrt{K'}}.
\ee Hence, by \cite[Lemma 2.4]{ZT}, we have

\be\label{KM-6}|P[\varphi_{0}]_{z}(z)|\geq\frac{1+K^{\ast}}{2K^{\ast}}\left(\frac{2}{\pi}-|P[\varphi_{0}](0)|\right),
\ee which, together with (\ref{KM-1}), yields that

\beq\label{KM-7}
\lambda(D_{f}(z))&\geq&\frac{2}{K+1}|f_{z}(z)|-\frac{\sqrt{K'}}{K+1}\\
\nonumber
&\geq&\frac{2}{K+1}\big(|P[\varphi_{0}]_{z}(z)|-\mathcal{H}(\varphi_{1},\cdots,\varphi_{n})\big)-\frac{\sqrt{K'}}{K+1}\\
\nonumber
&\geq&\frac{(1+K^{\ast})}{K^{\ast}(1+K)}\left(\frac{2}{\pi}-|P[\varphi_{0}](0)|\right)-\frac{2\mathcal{H}(\varphi_{1},\cdots,\varphi_{n})+\sqrt{K'}}{K+1}.
\eeq

\bcl\label{C-4.1}
$$\frac{(1+K^{\ast})}{K^{\ast}(1+K)}\left(\frac{2}{\pi}-|P[\varphi_{0}](0)|\right)-\frac{2\mathcal{H}(\varphi_{1},\cdots,\varphi_{n})+\sqrt{K'}}{K+1}>0.$$
\ecl Now we prove this Claim. Let
$\mathcal{B}=2/\pi-|P[\varphi_{0}](0)|$. By the assumptions, we have

\be\label{KM-8}2K\mathcal{H}(\varphi_{1},\cdots,\varphi_{n})+\sqrt{K'}<\mathcal{B}.\ee

It follows from (\ref{KM-5.1}) and (\ref{KM-8}) that

\beqq\label{KM-9}
\frac{K^{\ast}+1}{K^{\ast}(K+1)}\mathcal{B}&=&\frac{\mathcal{B}^{2}}{K\mathcal{B}+2K\mathcal{H}(\varphi_{1},\cdots,\varphi_{n})+\sqrt{K'}}
>\frac{\mathcal{B}^{2}}{K\mathcal{B}+\mathcal{B}}\\ \nonumber
&=&\frac{\mathcal{B}}{K+1}>\frac{2K\mathcal{H}(\varphi_{1},\cdots,\varphi_{n})+\sqrt{K'}}{K+1}\\
\nonumber
&\geq&\frac{2\mathcal{H}(\varphi_{1},\cdots,\varphi_{n})+\sqrt{K'}}{K+1},
\eeqq which implies that the Claim \ref{C-4.1} is true. Since for
all $z_{1},$ $z_{2}\in\mathbb{D}$,

\beqq|f(z_{1})-f(z_{2})|&\geq&
\int_{[z_{1},z_{2}]}\lambda(D_{f}(z))|dz|\\ &\geq&
\left(\frac{(1+K^{\ast})}{K^{\ast}(1+K)}\mathcal{B}-\frac{2\mathcal{H}(\varphi_{1},\cdots,\varphi_{n})+\sqrt{K'}}{K+1}\right)|z_{1}-z_{2}|,
\eeqq
 we conclude that $f$ is also co-Lipschitz continuous.

\bst\label{bst-4.2} The Lipschitz continuity of $f$.\est

The Lipschitz continuity of $f$ easy follows from (\ref{eq-x1.1}),
(\ref{eq-o1.1}) and Theorem \Ref{KM-P}.

Next, we prove part $(b)$.

\bst\label{bst-4.3} The asymptotically sharp Lipschitz inequality of
$f$. \est

 Since $P[\varphi_{0}]$ is a
$K^{\ast}$-quasiconformal mapping of $\mathbb{D}$ onto itself with
$P[\varphi_{0}](0)=0$, by \cite[Theorem 3.3]{PS}, we see that, for
all $z_{1},z_{2}\in\mathbb{D}$,

\be\label{KM-11}|P[\varphi_{0}](z_{1})-P[\varphi_{0}](z_{2})|\leq
(K^{\ast})^{3K^{\ast}+1}2^{5(K^{\ast}-1/K^{\ast})/2}|z_{1}-z_{2}|.
\ee By Lemmas \ref{lem-1} and \ref{lem-2},  we obtain that, for all
$z_{1},z_{2}\in\mathbb{D}$,

\be\label{KM-12}|G_{1}[\varphi_{1}](z_{1})-G_{1}[\varphi_{1}](z_{2})|\leq\frac{2}{3}\|\varphi_{1}\|_{\infty}|z_{1}-z_{2}|\ee
and
\be\label{KM-13}|G_{k}[\varphi_{k}](z_{1})-G_{k}[\varphi_{k}](z_{2})|\leq\frac{2}{15}\left(\frac{3}{16}\right)^{k-1}\|\varphi_{k}\|_{\infty}|z_{1}-z_{2}|,\ee
where $k\in\{2,\ldots,n\}$. It follows from (\ref{eq-ch-3.0}),
(\ref{KM-11}), (\ref{KM-12}) and (\ref{KM-13}) that, for all
$z_{1},z_{2}\in\mathbb{D}$,

\beqq|f(z_{1})-f(z_{2})|&\leq&|P[\varphi_{0}](z_{1})-P[\varphi_{0}](z_{2})|+\sum_{k=1}^{n}|G_{k}[\varphi_{k}](z_{1})-G_{k}[\varphi_{k}](z_{2})|\\
&\leq&\big(M_{1}(K,
K')+N_{1}(\varphi_{1},\ldots,\varphi_{n})\big)|z_{1}-z_{2}|,\eeqq
where $M_{1}(K,
K')=(K^{\ast})^{3K^{\ast}+1}2^{5(K^{\ast}-1/K^{\ast})/2}$ and
$$N_{1}(\varphi_{1},\ldots,\varphi_{n})=\frac{2}{3}\|\varphi_{1}\|_{\infty}+
\sum_{k=2}^{n}\frac{2}{15}\left(\frac{3}{16}\right)^{k-1}\|\varphi_{k}\|_{\infty}.$$
It is easy to know that

$$\lim_{K\rightarrow1^{+},K'\rightarrow0^{+}} M_{1}(K,
K')=1~\mbox{and}~\lim_{\|\varphi_{1}\|_{\infty}\rightarrow0^{+},\ldots,\|\varphi_{n}\|_{\infty}\rightarrow0^{+}}N_{1}(\varphi_{1},\ldots,\varphi_{n})\big)=0.$$

\bst\label{bst-4.4} The asymptotically sharp  co-Lipschitz
 inequality of $f$.\est
 Let $$M_{2}(K,
K')=\frac{(1+K^{\ast})}{K^{\ast}(1+K)}\max\left\{\frac{2}{\pi},L_{K^{\ast}}\right\}-\frac{\sqrt{K'}}{K+1}$$
and
$$N_{2}(\varphi_{1},\ldots,\varphi_{n})=\frac{2\mathcal{H}(\varphi_{1},\cdots,\varphi_{n})}{K+1},$$
where $K^{\ast}$ is defined in (\ref{KM-5.1}) and $L_{K^{\ast}}$ is
a positive constant satisfying
$$\lim_{K\rightarrow1^{+},K'\rightarrow0^{+}}L_{K^{\ast}}=1.$$ Then $$\lim_{K\rightarrow1^{+},K'\rightarrow0^{+}} M_{2}(K,
K')=1~\mbox{and}~\lim_{\|\varphi_{1}\|_{\infty}\rightarrow0^{+},\ldots,\|\varphi_{n}\|_{\infty}\rightarrow0^{+}}N_{2}(\varphi_{1},\ldots,\varphi_{n})\big)=0.$$
It follows from the Claim \ref{C-4.1} that

\be\label{KM-13.1}M_{2}(K,
K')-N_{2}(\varphi_{1},\ldots,\varphi_{n})>0.\ee Then, by
  Theorem \Ref{PS}, we have

\be\label{KM-15}|P[\varphi_{0}]_{z}(z)|\geq\frac{K^{\ast}+1}{2K^{\ast}}\max\left\{\frac{2}{\pi},L_{K^{\ast}}\right\},~z\in\mathbb{D},\ee
which, together with (\ref{KM-1}) and (\ref{KM-13.1}), yields that,
for all $z_{1},z_{2}\in\mathbb{D}$,

\beqq\label{KM-17} |f(z_{1})-f(z_{2})|&\geq&
\int_{[z_{1},z_{2}]}\lambda(D_{f}(z))|dz|\geq\int_{[z_{1},z_{2}]}\left(\frac{2}{K+1}|f_{z}(z)|-\frac{\sqrt{K'}}{K+1}\right)|dz|\\
\nonumber
&\geq&\int_{[z_{1},z_{2}]}\left(\frac{2}{K+1}\big(|P[\varphi_{0}]_{z}(z)|-\mathcal{H}(\varphi_{1},\cdots,\varphi_{n})\big)-\frac{\sqrt{K'}}{K+1}\right)|dz|\\
\nonumber &\geq&\big(M_{2}(K,
K')-N_{2}(\varphi_{1},\ldots,\varphi_{n})\big)|z_{1}-z_{2}|.\eeqq
Therefore, $f$ is co-Lipschitz continuous in $\mathbb{D}$. The proof
of this theorem is complete.\qed

\section{The proof of Theorem \ref{thm-1.1}}\label{sec-4}
The purpose of this section is to prove Theorem \ref{thm-1.1}. The
proof consists of three steps. In the first step, the Lipschitz continuity of the mappings $f$ is proved, the co-Lipschitz continuity of $f$ is demonstrated in the second step,
and in the third step, the Lipschitz and co-Lipschitz continuity coefficients obtained in the first two steps are shown to have bounds
with the forms as required in Theorem \ref{thm-1.1}.

\bst\label{bst-1} The asymptotically sharp Lipschitz inequality of
$f$.
 \est

We start the discussions of this step with the following claim. \bcl
The limits $$\lim_{z\rightarrow
\xi\in\mathbb{T},z\in\mathbb{D}}D_{f}(z)\;\; \mbox{and}\;\;
\lim_{z\rightarrow \xi\in\mathbb{T},z\in\mathbb{D}}J_{f}(z)$$ exist
almost everywhere in $\mathbb{T}$. \ecl

We are going to verify the existence of  these two limits by
applying Theorem \Ref{KM-P} and Lemma \ref{lem-main}. For this, we
need to get an upper bound of  $|\Delta f|$ as stated
in \eqref{eq-x1.1} and \eqref{eq-o1.1} below, and we will divide it into two cases to
estimate.

\noindent $\mathbf{Case~ 1.}$  $n=2.$

By the formula (1.3) in \cite{K2} (see also \cite[pp. 118-120]{Ho}),
we have that for $z\in \ID$, \beqq\label{eq-yy2}\Delta
f(z)=P[\varphi_{1}](z)-\int_{\mathbb{D}}G(z,\zeta)\varphi_{2}(\zeta)d\sigma(\zeta).
\eeqq It follows from Lemma \ref{lem-0.1} (\ref{kk-1}) that
\be\label{eq-x1.1} |\Delta f(z)|\leq
|P[\varphi_{1}](z)|+\|\varphi_{2}\|_{\infty}\int_{\mathbb{D}}|G(z,\zeta)|d\sigma(\zeta)
\leq \|\varphi_{1}\|_{\infty}+\frac{\|\varphi_{2}\|_{\infty}}{4}.\ee

\noindent $\mathbf{Case~ 2.}$  $n\geq3.$

Since $$\Delta^{n-1}(\Delta f)=\varphi_{n}~\mbox{in}~\mathbb{D},$$
and $$\Delta^{n-2}(\Delta f)=\varphi_{n-1},~\ldots,~\Delta f
=\varphi_{1}$$ on $\mathbb{T}$, by (\ref{eq-ch-3.0}), we see that,
for $z\in \mathbb{D}$,
\begin{eqnarray*}
\Delta
f(z)=P[\varphi_{1}](z)+\sum_{j=1}^{n-1}(-1)^{j}G_{j}[\varphi_{j+1}](z),
\end{eqnarray*} where
\begin{eqnarray*}
G_{k}[\varphi_{k+1}](z)&=&
\int_{\mathbb{D}}\cdots\int_{\mathbb{D}}G(z,\xi_{1})\cdots
G(\xi_{k-1},\xi_{k})\\ \nonumber
&&\times\left(\int_{0}^{2\pi}P(\xi_{k},e^{it})\varphi_{k+1}(e^{it})dt\right)d\sigma(\xi_{k})\cdots
d\sigma(\xi_{1})
\end{eqnarray*}for $k\in\{1,\ldots,n-2\}$, and

\begin{eqnarray*}
 G_{n-1}[\varphi_{n}](z)&=&
\int_{\mathbb{D}}\cdots\int_{\mathbb{D}}G(z,\zeta_{1})\cdots G(\zeta_{n-3},\zeta_{n-2})\\
\nonumber
&&\times\left(\int_{\mathbb{D}}G(\zeta_{n-2},\zeta_{n-1})\varphi_{n}(\zeta_{n-1})d\sigma(\zeta_{n-1})\right)d\sigma(\zeta_{n-2})\cdots
d\sigma(\zeta_{1}).
\end{eqnarray*}

By Lemmas \ref{lem-0.1} and \Ref{Lem-es3}, for $z\in\mathbb{D}$, we
obtain that

\begin{eqnarray*}
|G_{k}[\varphi_{k+1}](z)|&\leq&\|\varphi_{k+1}\|_{\infty}
\int_{\mathbb{D}}\cdots\int_{\mathbb{D}}|G(z,\xi_{1})|\cdots
|G(\xi_{k-1},\xi_{k})|d\sigma(\xi_{k})\cdots d\sigma(\xi_{1})\\
&\leq&\frac{\|\varphi_{k+1}\|_{\infty}}{4}\left(\frac{3}{16}\right)^{k-1}(1-|z|^{2})\\
&\leq&\frac{\|\varphi_{k+1}\|_{\infty}}{4}\left(\frac{3}{16}\right)^{k-1}
\end{eqnarray*}
for $k\in\{1,\ldots,n-2\}$, and

\begin{eqnarray*}
 |G_{n-1}[\varphi_{n}](z)|&=&\|\varphi_{n}\|_{\infty}
\int_{\mathbb{D}}\cdots\int_{\mathbb{D}}|G(z,\zeta_{1})|\cdots
|G(\zeta_{n-3},\zeta_{n-2})|d\sigma(\zeta_{n-2})\cdots
d\sigma(\zeta_{1})\\
&\leq&\frac{\|\varphi_{n}\|_{\infty}}{4}\left(\frac{3}{16}\right)^{n-2}(1-|z|^{2})\\
&\leq&\frac{\|\varphi_{n}\|_{\infty}}{4}\left(\frac{3}{16}\right)^{n-2},
\end{eqnarray*}
which give that


\beq\label{eq-o1.1} |\Delta
f(z)|&=&|P[\varphi_{1}](z)|+\sum_{j=1}^{n-1}|G_{j}[\varphi_{j+1}](z)|\\
\nonumber
&\leq&\|\varphi_{1}\|_{\infty}+\sum_{j=1}^{n-1}\frac{\|\varphi_{j+1}\|_{\infty}}{4}\left(\frac{3}{16}\right)^{j-1}<\infty.
\eeq Since $f$ is a $K$-quasiconformal self-mapping of $\mathbb{D}$,
we see that $f$ can  be extended to the homeomorphism of
$\overline{\mathbb{D}}$ onto itself. Now, the existence of the
limits
$$D_{f}(\xi)=\lim_{z\rightarrow
\xi\in\mathbb{T},z\in\mathbb{D}}D_{f}(z)~\mbox{and}~J_{f}(\xi)=\lim_{z\rightarrow
\xi\in\mathbb{T},z\in\mathbb{D}}J_{f}(z)$$ almost everywhere in
$\mathbb{T}$ follows from (\ref{eq-o1.1}), Theorem \Ref{KM-P}  and Lemma
\ref{lem-main}.

For convenience, in the following, let
$$ C_{2}(K,\varphi_{1},\cdots,\varphi_{n})=\sup_{z\in\mathbb{D}}\|D_{f}(z)\|.$$

Since for almost all $z_1$ and $z_2\in \ID$, \be\label{sun-7}
|f(z_1)-f(z_2)=\Big|\int_{[z_1,z_2]}f_zdz+f_{\overline{z}}d\overline{z}
\Big|\leq C_{2}(K,\varphi_{1},\cdots,\varphi_{n}) |z_1-z_2|, \ee we
see that, to prove the Lipschitz continuity of $f$ and investigate the behavior of the Lipschitz coefficient, it suffices to
estimate the quantity $C_{2}(K,\varphi_{1},\cdots,\varphi_{n})$. To
reach this goal, we first show that the quantity
$C_{2}(K,\varphi_{1},\cdots,\varphi_{n})$ satisfies an inequality
which is stated in the following claim. \bcl\label{eq-sh-18}
$C_{2}(K,\varphi_{1},\cdots,\varphi_{n})\leq
\big(C_{2}(K,\varphi_{1},\cdots,\varphi_{n})\big)^{1-\frac{1}{K}}\mu_1+\mu_2,$
 where $$\mu_1=\frac{K(Q(K))^{\frac{1}{K}+1}}{2\pi}\int_{0}^{2\pi}|1-e^{it}|^{-1+\frac{1}{K^{2}}}
dt<\infty,$$ $Q(K)$ is from Theorem \Ref{Mori}, $\mu_2=\mu_3+\mu_4,$
 $$\mu_3=\frac{K\|\varphi_{1}\|_{\infty}}{2}
 +K\sum_{k=2}^{n}\frac{\|\varphi_{k}\|_{\infty}}{16}\left(\frac{3}{16}\right)^{k-2},$$
 and
 $$\mu_4=\frac{7}{6}\|\varphi_{1}\|_{\infty}+\sum_{k=2}^{n}\frac{47\|\varphi_{k}\|_{\infty}}{240}\left(\frac{3}{16}\right)^{k-2}.$$
\ecl

To prove the claim, we need the following preparation. Firstly, we
prove that $\mu_1<\infty$, and for almost every $\theta\in [0,2\pi]$, \be\label{sun-3}
\|D_{f}(e^{i\theta})\|\leq\frac{K}{2\pi}\int_{0}^{2\pi}\frac{|f(e^{it})-f(e^{i\theta})|^{2}}{|e^{it}-e^{i\theta}|^{2}}dt+\mu_3.\ee
By
\cite[Lemma 1.6]{K3}, we know that
$$\int_{0}^{2\pi}|e^{it}-e^{i\theta_{\epsilon}}|^{-1+\frac{1}{K^{2}}}
dt<\infty,$$ which shows $\mu_1<\infty$. Next, we prove (\ref{sun-3}).
 For $\theta\in [0,2\pi]$, let
$$\varphi_{0}(e^{i\theta})=f(e^{i\theta})=e^{i\gamma(\theta)}.$$
Then, by (\ref{eq-d1}), we see that $f(e^{i\theta})$ is absolutely
continuous. It follows that
$$i\gamma'(\theta)e^{i\gamma(\theta)}=\frac{d}{d\theta}f(e^{i\theta})=\lim_{r\rightarrow1^{-}}\frac{\partial}{\partial\theta}f(re^{i\theta})=\lim_{r\rightarrow1^{-}}
\big[ir\left(f_{z}(re^{i\theta})e^{i\theta}-f_{\overline{z}}(re^{i\theta})e^{-i\theta}\right)\big],$$
which implies
\be\label{eq-x1.2}\frac{1}{K}\|D_{f}(e^{i\theta})\|\leq
\lim_{r\rightarrow1^{-}}\lambda(D_{f}(re^{i\theta}))\leq\gamma'(\theta)\leq\lim_{r\rightarrow1^{-}}\|D_{f}(re^{i\theta})\|=\|D_{f}(e^{i\theta})\|\ee
almost everywhere in $[0,2\pi],$ where $r\in[0,1)$.

Since the existence of the two limits
$$D_{f}(e^{i\theta})=\lim_{z\rightarrow
e^{i\theta},z\in\mathbb{D}}D_{f}(z)~\mbox{and}~J_{f}(e^{i\theta})=\lim_{z\rightarrow
e^{i\theta},z\in\mathbb{D}}J_{f}(z)$$ almost everywhere in
$[0,2\pi]$ guarantees that
$$\|D_{f}(e^{i\theta})\|^{2}\leq K J_{f}(e^{i\theta}),$$
we deduce from (\ref{eq-sh-1}) and (\ref{eq-x1.2}) that
\beq\label{eq-sh-11}\nonumber
 \|D_{f}(e^{i\theta})\|^{2}\leq K
\|D_{f}(e^{i\theta})\|\Bigg\{\frac{1}{2\pi}\int_{0}^{2\pi}\frac{|f(e^{it})-f(e^{i\theta})|^{2}}{|e^{it}-e^{i\theta}|^{2}}dt+\frac{\mu_3}{K}\Bigg\},\eeq
from which the inequality \eqref{sun-3} follows.

Secondly, we show that for any $\epsilon>0$, there exists
$\theta_{\epsilon}\in [0,2\pi]$ such that \be\label{sun-5}
C_{2}(K,\varphi_{1},\cdots,\varphi_{n})\leq(1+\epsilon)\|D_{f}(e^{i\theta_{\epsilon}})\|+\mu_4.\ee

For the proof, let $t\in[0,2\pi]$, and let
$$H_{t}(z)=\frac{\partial}{\partial
z}P[\varphi_{0}](z)+e^{it}\overline{\frac{\partial}{\partial
\overline{z}}P[\varphi_{0}](z)}$$ in $\mathbb{D}$.


Since $P[\varphi_{0}]=f-\sum_{k=1}^{n}(-1)^{k}G_{k}[\varphi_{k}]$ is
harmonic, we see that $H_{t}$ is analytic in $\mathbb{D}$, and thus,
$$|H_{t}(z)|\leq\mbox{esssup}_{\theta\in[0,2\pi]}|H_{t}(e^{i\theta})|\leq\mbox{esssup}_{\theta\in[0,2\pi]}\|D_{P[\varphi_{0}]}(e^{i\theta})\|.$$
Then, the facts
$$\|D_{P[\varphi_{0}]}(z)\|=\max_{t\in[0,2\pi]}|H_{t}(z)|\leq\mbox{esssup}_{\theta\in[0,2\pi]}\|D_{P[\varphi_{0}]}(e^{i\theta})\|$$
and
$$\|D_{P[\varphi_{0}]}(z)\|=\left|\frac{\partial f}{\partial z}-\sum_{k=1}^{n}(-1)^{k}\frac{\partial}{\partial z}G_{k}[\varphi_{k}]\right|+
\left|\frac{\partial f}{\partial
\overline{z}}-\sum_{k=1}^{n}(-1)^{k}\frac{\partial}{\partial
\overline{z}}G_{k}[\varphi_{k}]\right|
$$
 ensure

\begin{eqnarray*}
\|D_{P[\varphi_{0}]}(z)\|&\leq&
\mbox{esssup}_{\theta\in[0,2\pi]}\|D_{f}(e^{i\theta})\|+
\sum_{k=1}^{n}\mbox{esssup}_{\theta\in[0,2\pi]}\|D_{G_{k}[\varphi_{k}]}(e^{i\theta})\|,
\end{eqnarray*} which, together with Lemmas \ref{lem-1} and
\ref{lem-2}, guarantees that for all $z\in\mathbb{D},$
$$
\|D_{f}(z)\|\leq\|D_{P[\varphi_{0}]}(z)\|+\sum_{k=1}^{n}\|D_{G_{k}[\varphi_{k}]}(z)\|\leq
\mbox{esssup}_{\theta\in[0,2\pi]}\|D_{f}(e^{i\theta})\|+\mu_4,
$$ from which the inequality \eqref{sun-5}  follows.

Let
$$\nu=\frac{1}{2\pi}\int_{0}^{2\pi}\frac{|f(e^{it})-f(e^{i\theta_{\epsilon}})|^{2}}{|e^{it}-e^{i\theta_{\epsilon}}|^{2}}dt.$$
Finally, we need the following estimate of $\nu$:
\be\label{sun-6}\nu\leq
\frac{\big(C_{2}(K,\varphi_{1},\cdots,\varphi_{n})\big)^{1-\frac{1}{K}}(Q(K))^{\frac{1}{K}+1}}{2\pi}\int_{0}^{2\pi}|e^{it}-e^{i\theta_{\epsilon}}|^{-1+\frac{1}{K^{2}}}
dt. \ee

Since it follows from \eqref{sun-7} that for almost all
$\theta_{1}$, $\theta_{2}\in [0,2\pi],$
\be\label{sun-8}|f(e^{i\theta_{1}})-f(e^{i\theta_{2}})|\leq
C_{2}(K,\varphi_{1},\cdots,\varphi_{n})\left|e^{i\theta_{1}}-e^{i\theta_{2}}\right|,\ee
we infer that
$$ \nu \leq
\frac{\big(C_{2}(K,\varphi_{1},\cdots,\varphi_{n})\big)^{1-\frac{1}{K}}}{2\pi}\int_{0}^{2\pi}|e^{it}-e^{i\theta_{\epsilon}}|^{-1+\frac{1}{K^{2}}}
\frac{|f(e^{it})-f(e^{i\theta_{\epsilon}})|^{1+\frac{1}{K}}}{|e^{it}-e^{i\theta_{\epsilon}}|^{\frac{1}{K}+\frac{1}{K^{2}}}}dt,
$$ from which, together with   Theorem \Ref{Mori}, the inequality \eqref{sun-6} follows.\medskip

Now, we are ready to finish the proof of the claim. It follows from
 \eqref{sun-5} that
$$C_{2}(K,\varphi_{1},\cdots,\varphi_{n})\leq(1+\epsilon)\|D_{f}(e^{i\theta_{\epsilon}})\|+\mu_4,$$
and so, \eqref{sun-3} and \eqref{sun-6} give \be\label{eq-sh-17}
C_{2}(K,\varphi_{1},\cdots,\varphi_{n}) \leq
\big(C_{2}(K,\varphi_{1},\cdots,\varphi_{n})\big)^{1-\frac{1}{K}}\mu_1(1+\epsilon)+\mu_3(1+\epsilon)+\mu_4.
 \ee

By letting $\epsilon\rightarrow0^{+}$, we get from (\ref{eq-sh-17})
that
$$C_{2}(K,\varphi_{1},\cdots,\varphi_{n})\leq
\big(C_{2}(K,\varphi_{1},\cdots,\varphi_{n})\big)^{1-\frac{1}{K}}\mu_1+\mu_2,$$
as required.\medskip

The following is a lower bound for
$C_{2}(K,\varphi_{1},\cdots,\varphi_{n})$. \bcl\label{sun-1}
$C_{2}(K,\varphi_{1},\cdots,\varphi_{n})\geq 1$.\ecl

 Since
$$\int_{0}^{2\pi}\gamma'(\theta)d\theta=\gamma(2\pi)-\gamma(0)=2\pi,$$
we conclude that
$$\mbox{esssup}_{\theta\in[0,2\pi]}\lim_{t\rightarrow\theta}\left|\frac{f(e^{i\theta})-f(e^{it})}{e^{i\theta}-e^{it}}\right|=\mbox{esssup}_{\theta\in[0,2\pi]}\gamma'(\theta)\geq1.$$
Then, it follows from (\ref{sun-8}) and the following fact
$$\mbox{esssup}_{\theta\in[0,2\pi]}\lim_{t\rightarrow\theta}\left|\frac{f(e^{i\theta})-f(e^{it})}{e^{i\theta}-e^{it}}\right| \leq \mbox{esssup}_{0\leq\theta\neq
 t\leq 2\pi}\left|\frac{f(e^{i\theta})-f(e^{it})}{e^{i\theta}-e^{it}}\right|
 $$  that
$$C_{2}(K,\varphi_{1},\cdots,\varphi_{n})\geq 1.$$ Hence, the claim is true. \medskip

An upper bound of $C_{2}(K,\varphi_{1},\cdots,\varphi_{n})$ is
established in the following claim. \bcl\label{sun-2} If
$\frac{(K-1)}{K}\mu_1<1$, then
$$C_{2}(K,\varphi_{1},\cdots,\varphi_{n})\leq \mu_5,$$ where
$\mu_5=\frac{\frac{1}{K}\mu_1+\mu_2}{1-\mu_1\left(1-\frac{1}{K}\right)}.$
\ecl
 The proof of this claim easily follows from \cite[Lemma 2.9]{K2}.\medskip

Now, we are ready to finish the discussions in this step. By Claims \ref{eq-sh-18} and \ref{sun-1}, we obtain
$$1\leq C_{2}(K,\varphi_{1},\cdots,\varphi_{n})\leq \mu_6,$$ where $\mu_6=(\mu_1+\mu_2)^{K}.$

By letting
$$C_{3}=\begin{cases}
\displaystyle \mu_6, &\mbox{ if }\, \frac{(K-1)}{K}\mu_1\geq1,\\
\displaystyle \min\{\mu_5, \mu_6\}, &\mbox{ if }\,
\frac{(K-1)}{K}\mu_1<1,
\end{cases}$$ we infer that
\be\label{eq-cw1} 1<C_{2}(K,\varphi_{1},\cdots,\varphi_{n})\leq
C_{3}.\ee Then, the Lipschtz continuity of $f$  follows from these
estimates of $C_{2}(K,\varphi_{1},\cdots,\varphi_{n})$.

\bst\label{bst-2} The asymptotically sharp co-Lipschitz inequality
of $f$. \est We begin the discussions of this step with some
preparation which consists of the following two claims.

\bcl\label{eq-sh-28} $\lambda(D_{P[\varphi_{0}]}(e^{i\theta}))\geq
\frac{\mu_7}{K^{2}}-\left(1+\frac{1}{K^{2}}\right)\mu_8$ almost everywhere on $\theta\in[0,2\pi]$,
 \noindent where
 \be\label{thur-2} \mu_7=\max\{\mu'_7,~\mu''_7\},\;\;
\mu'_7=(Q(K))^{-2K}\frac{1}{2\pi}\int_{0}^{2\pi}|e^{it}-e^{i\theta}|^{2K-2}dt,
\ee
$$\mu''_7=\frac{1}{2}-\sum_{k=1}^{n}\frac{\|\varphi_{k}\|_{\infty}}{8}\left(\frac{3}{16}\right)^{k-1},$$ and
\be\label{thur-1} \mu_8=\frac{\|\varphi_{1}\|_{\infty}}{2}
+\sum_{k=2}^{n}\frac{\|\varphi_{k}\|_{\infty}}{16}\left(\frac{3}{16}\right)^{k-2}.
\ee \ecl

By (\ref{eq-x1.2}), we have

 $$\frac{\gamma'(\theta)}{K}\leq\frac{\|D_{f}(e^{i\theta})\|}{K}\leq\lambda(D_{f}(e^{i\theta}))
 \leq\lambda(D_{P[\varphi_{0}]}(e^{i\theta}))+\sum_{k=1}^{n}\|D_{G_{k}[\varphi_{k}]}(e^{i\theta})\|,$$
which, together with Lemmas \ref{lem-1} and \ref{lem-2}, implies

\be\label{eerr-1}\lambda(D_{P[\varphi_{0}]}(e^{i\theta}))\geq
\frac{\gamma'(\theta)}{K}
-\sum_{k=1}^{n}\|D_{G_{k}[\varphi_{k}]}(e^{i\theta})\|\geq\frac{\gamma'(\theta)}{K}
-\mu_8.\ee Then, we know from (\ref{eerr-1})  that, to prove the
claim, it suffices to show that \be\label{mon-2}
K\gamma'(\theta)\geq \mu_7.\ee

Again, it follows from (\ref{eq-x1.2}) that
$$\frac{J_{f}(e^{i\theta})}{\gamma'(\theta)}\leq\frac{J_{f}(e^{i\theta})}{\lambda(D_{f}(e^{i\theta}))}
\leq K\lambda(D_{f}(e^{i\theta}))\leq K\gamma'(\theta),$$ and thus,
(\ref{eq-sh-1.1}) gives
$$K\gamma'(\theta)\geq\frac{1}{2\pi}\int_{0}^{2\pi}\frac{|f(e^{it})-f(e^{i\theta})|^{2}}{|e^{it}-e^{i\theta}|^{2}}dt
-\mu_8. $$ This implies that, to prove \eqref{mon-2}, we only need
to verify the validity of the following inequality:
\be\label{eq-sh-25}\frac{1}{2\pi}\int_{0}^{2\pi}\frac{|f(e^{it})-f(e^{i\theta})|^{2}}{|e^{it}-e^{i\theta}|^{2}}dt
\geq \mu_7. \ee

We now prove this inequality. On the one hand, since $f^{-1}$ is  a
$K$-quasiconformal mapping, it follows from Theorem \Ref{Mori} that
for any $z_{1},z_{2}\in\mathbb{D}$, \beqq\label{eq-sh-26}
(Q(K))^{-K}|z_{1}-z_{2}|^{K}\leq|f(z_{1})-f(z_{2})|,\eeqq
which implies \beq\label{eq-sh-27}
\frac{1}{2\pi}\int_{0}^{2\pi}\frac{|f(e^{it})-f(e^{i\theta})|^{2}}{|e^{it}-e^{i\theta}|^{2}}dt\geq
\mu'_7.\eeq

On the other hand, since $f(0)=0$, we see from
\begin{eqnarray*}
 |G_{k}[\varphi_{k}](0)|&=&\|\varphi_{k}\|_{\infty}
\int_{\mathbb{D}}\cdots\int_{\mathbb{D}}|G(0,\zeta_{1})|\cdots
|G(\xi_{k-1},\xi_{k})|d\sigma(\xi_{k})\cdots
d\sigma(\xi_{1})\\
&\leq&\frac{\|\varphi_{k}\|_{\infty}}{4}\left(\frac{3}{16}\right)^{k-1}
\end{eqnarray*}
and

\begin{eqnarray*}
 |G_{n}[\varphi_{n}](0)|&=&\frac{\|\varphi_{n}\|_{\infty}}{4}
\int_{\mathbb{D}}\cdots\int_{\mathbb{D}}|G(0,\zeta_{1})|\cdots
|G(\zeta_{n-2},\zeta_{n-1})|\\
&&\times(1-|\zeta_{n-1}|^{2})d\sigma(\zeta_{n-1})\cdots
d\sigma(\zeta_{1})\\
&\leq&\frac{\|\varphi_{n}\|_{\infty}}{4}\left(\frac{3}{16}\right)^{n-1}
\end{eqnarray*}
that \be\label{ccww-1} |P[\varphi_{0}](0)|
\leq\sum_{k=1}^{n}|G_{k}[\varphi_{k}](0)|\leq
\sum_{k=1}^{n}\frac{\|\varphi_{k}\|_{\infty}}{4}\left(\frac{3}{16}\right)^{k-1}.
\ee Then, we infer from (\ref{ccww-1}) and the following fact:
$$\frac{1}{2\pi}\int_{0}^{2\pi}\frac{|f(e^{it})-f(e^{i\theta})|^{2}}{|e^{it}-e^{i\theta}|^{2}}dt\geq
\frac{1}{4\pi}\int_{0}^{2\pi}\left[1-\mbox{Re}\big(f(e^{it})\overline{f(e^{i\theta})}\big)\right]dt\geq\frac{1-|P[\varphi_{0}](0)|}{2}$$
that \beq\label{ccww-2}
\frac{1}{2\pi}\int_{0}^{2\pi}\frac{|f(e^{it})-f(e^{i\theta})|^{2}}{|e^{it}-e^{i\theta}|^{2}}dt
\geq \mu''_7.\eeq

Obviously, the inequality \eqref{eq-sh-25} follows from
\eqref{eq-sh-27} and \eqref{ccww-2}, and so, the claim is proved.
\bcl\label{ccww-6} For $z\in\mathbb{D},$
$\lambda(D_{P[\varphi_{0}]}(z))\geq
\frac{\mu_7}{K^{2}}-\left(1+\frac{1}{K^{2}}\right)\mu_8.$ \ecl


By the Choquet-Rad\'o-Kneser theorem (see \cite{AJK,Cho}), we see
that $P[\varphi_{0}]$ is a  sense-preserving  harmonic
diffeomorphism of $\mathbb{D}$ onto itself. It follows from  Lewy's
theorem (cf. \cite{Lw}) and \cite[Inequality (17)]{He} that

\be\label{dwc-I}\inf_{z\in\mathbb{D}}\left|\frac{\partial}{\partial
z}P[\varphi_{0}](z)\right|>0.\ee

Hence, for $z\in\mathbb{D}$, we can let
$$p_{1}(z)=\frac{\overline{\frac{\partial}{\partial
\overline{z}}P[\varphi_{0}](z)}}{\frac{\partial}{\partial
z}P[\varphi_{0}](z)}~\mbox{and}~p_{2}(z)=\left(\frac{\mu_7}{K^{2}}-\mu_8\right)\frac{1}{\frac{\partial}{\partial
z}P[\varphi_{0}](z)},$$ and  let
$$q_{\vartheta}(z)=p_{1}(z)+e^{i\vartheta}p_{2}(z),$$ where  $\vartheta\in[0,2\pi]$.
Since $P[\varphi_{0}]$ is a  sense-preserving  harmonic
diffeomorphism of $\mathbb{D}$, by (\ref{dwc-I}), we see that

\be\label{dwc-II}\sup_{z\in\mathbb{D}}|q_{\vartheta}(z)|<+\infty.\ee

By Claim \ref{eq-sh-28}, we have

\be\label{dwc-III}|q_{\vartheta}(e^{i\theta})|\leq|p_{1}(e^{i\theta})|+|p_{2}(e^{i\theta})|=
\frac{\left|\frac{\partial}{\partial
\overline{z}}P[\varphi_{0}](e^{i\theta})\right|+\frac{\mu_7}{K^{2}}-\left(1+\frac{1}{K^{2}}\right)\mu_8}
{\left|\frac{\partial}{\partial
z}P[\varphi_{0}](e^{i\theta})\right|}\leq1\ee almost everywhere on
$\theta\in[0,2\pi].$

Let $$E=\{\theta\in[0,2\pi]:~\lim_{z\rightarrow
e^{i\theta}}q_{\vartheta}(z)~\mbox{exists}\}.$$ Then the measure of
the set $[0,2\pi]\backslash E$ is zero. Hence, for $r\in[0,1)$, we
have

\begin{eqnarray*}
|q_{\vartheta}(rz)|&\leq&\int_{0}^{2\pi}P(z,e^{i\theta})|q_{\vartheta}(re^{i\theta})|d\theta\\
&\leq&\int_{E}P(z,e^{i\theta})|q_{\vartheta}(re^{i\theta})|d\theta+\int_{[0,2\pi]\backslash
E}P(z,e^{i\theta})|q_{\vartheta}(re^{i\theta})|d\theta
\end{eqnarray*}
which, together with (\ref{dwc-II}), (\ref{dwc-III}) and  the
Lebesgue Dominated Convergence Theorem, implies
\be\label{cckkww-7}|q_{\vartheta}(z)|\leq
\int_{E}P(z,e^{i\theta})|q_{\vartheta}(e^{i\theta})|d\theta
\leq1,\ee where $z\in\mathbb{D}$. It follows from (\ref{cckkww-7})
and the arbitrariness of $\vartheta\in[0,2\pi]$ that, for
$z\in\mathbb{D}$,
$$|p_{1}(z)|+|p_{2}(z)|\leq1,$$
from which the claim follows.

Now, we are ready to finish the proof of the co-Lipschitz continuity
of $f$. Since
$$\lambda(D_{f}(z))\geq\lambda(D_{P[\varphi_{0}]}(z))-\sum_{k=1}^{n}\|D_{G_{k}[\varphi_{k}]}(z)\|,$$
we see from Claim \ref{ccww-6}, Lemmas \ref{lem-1} and \ref{lem-2}
that \beq\label{mon-3}\lambda(D_{f}(z))\geq
C_{1}(K,\varphi_{1},\cdots,\varphi_{n}), \eeq where
\beq\label{ccww-10}
C_{1}(K,\varphi_{1},\cdots,\varphi_{n})&=&\frac{\mu_7}{K^{2}}-\left(1+\frac{1}{K^{2}}\right)\mu_8\\
\nonumber&&-\frac{2}{3}\|\varphi_{1}\|_{\infty}-
\sum_{k=2}^{n}\frac{2\|\varphi_{k}\|_{\infty}}{15}\left(\frac{3}{16}\right)^{k-2}.
\eeq And, we know from  \eqref{thur-2} and \eqref{thur-1} that
$C_{1}(K,\varphi_{1},\cdots,\varphi_{n})>0$ for small enough
 $\|\varphi_{k}\|_{\infty}$, where $k\in\{1,2,\ldots,n\}$. Since for all $z_{1},$
$z_{2}\in\mathbb{D}$,
$$|f(z_{1})-f(z_{2})|\geq \int_{[z_{1},z_{2}]}\lambda(D_{f}(z))|dz|\geq C_{1}(K,\varphi_{1},\cdots,\varphi_{n})|z_{1}-z_{2}|,$$
 we conclude that $f$ is co-Lipschitz continuous.

\bst\label{bst-3} Bounds of the Lipschitz continuity coefficients
$C_{1}(K,\varphi_{1},\cdots,\varphi_{n})$ and
$C_{2}(K,\varphi_{1},\cdots,\varphi_{n})$.\est

The discussions of this step consists of the following two claims.

\bcl There are constants $M_{3}(K)$ and
$N_{3}(K,\varphi_{1},\cdots,\varphi_{n})$ such that
\begin{enumerate} \item
$C_{2}(K,\varphi_{1},\cdots,\varphi_{n})\leq
M_{3}(K)+N_{3}(K,\varphi_{1},\cdots,\varphi_{n});$
\item
$\lim_{K\rightarrow1^{+}}M_{3}(K)=1,$ and
\item
$$\lim_{\|\varphi_{1}\|_{\infty}\rightarrow0^{+},\cdots,\|\varphi_{n}\|_{\infty}\rightarrow0^{+}}N_{3}(K,\varphi_{1},\cdots,\varphi_{n})=0.$$
\end{enumerate}\ecl
 From   (\ref{eq-cw1}), we see
that
$$ 1\leq C_{2}(K,\varphi_{1},\cdots,\varphi_{n})\leq C_3,$$
where
$$ C_3= \begin{cases}
\displaystyle (\mu_1+\mu_2)^{K}, &\mbox{ if }\, \frac{(K-1)}{K}\mu_1\geq 1,\\
\displaystyle \min\left\{(\mu_1+\mu_2)^{K},
\frac{\frac{1}{K}\mu_1+\mu_2}{1-\mu_1\left(1-\frac{1}{K}\right)}\right\},
&\mbox{ if }\, \frac{(K-1)}{K}\mu_1<1.
\end{cases}$$
Then, we have

$$C_3=\begin{cases}
\displaystyle M^{\ast}_{1}, &\mbox{ if }\, \frac{(K-1)}{K}\mu_1\geq 1,\\
\displaystyle \min\big\{M^{\ast}_{1}, M^{\ast}_{2}\big\}, &\mbox{ if
}\, \frac{(K-1)}{K}\mu_1< 1,
\end{cases}$$
where
$M^{\ast}_{1}=M_{2}^{'}(K)+N_{2}^{'}(K,\varphi_{1},\cdots,\varphi_{n})$,
$M^{\ast}_{2}=M_{2}^{''}(K)+N_{2}^{''}(K,\varphi_{1},\cdots,\varphi_{n})$,
$M_{2}^{'}(K)=\mu_1^{K}$,
$M_{2}^{''}(K)=\frac{\mu_1}{K-\mu_1(K-1)}$,
$N_{2}^{'}(K,\varphi_{1},\cdots,\varphi_{n})=(\mu_1+\mu_2)^{K}-\mu_1^{K}$,
and
$$N_{2}^{''}(K,\varphi_{1},\cdots,\varphi_{n})=\frac{\mu_2}{1-\mu_1\big(1-\frac{1}{K}\big)}.$$

Let
$$  M_{3}(K)=\begin{cases}
\displaystyle M_{2}^{'}(K), &\mbox{ if }\, \big(1-K^{-1}\big)\mu_1\geq1,\\
\displaystyle M_{2}^{''}(K), &\mbox{ if }\,
\big(1-K^{-1}\big)\mu_1<1 ~\mbox{and}~M_{1}^{\ast}\geq
M_{2}^{\ast},\\
\displaystyle M_{2}^{'}(K), &\mbox{ if }\, \big(1-K^{-1}\big)\mu_1<1
~\mbox{and}~M_{1}^{\ast}\leq M_{2}^{\ast}
\end{cases}$$
and
$$ N_{3}(K,\varphi_{1},\cdots,\varphi_{n})=\begin{cases}
\displaystyle N_{2}^{'}(K,\varphi_{1},\cdots,\varphi_{n}), &\mbox{ if }\, \big(1-K^{-1}\big)\mu_1\geq1,\\
\displaystyle N_{2}^{''}(K,\varphi_{1},\cdots,\varphi_{n}), &\mbox{
if }\, \big(1-K^{-1}\big)\mu_1<1 ~\mbox{and}~M_{1}^{\ast}\geq
M_{2}^{\ast},\\
\displaystyle N_{2}^{'}(K,\varphi_{1},\cdots,\varphi_{n}), &\mbox{
if }\, \big(1-K^{-1}\big)\mu_1<1 ~\mbox{and}~M_{1}^{\ast}\leq
M_{2}^{\ast}.
\end{cases}$$

It follows from the facts
$$\lim_{K\rightarrow1^{+}}M_{3}(K)=1\;\;\mbox{and}\;\; \lim_{\|\varphi_{1}\|_{\infty}\rightarrow0^{+},\cdots,\|\varphi_{n}\|_{\infty}\rightarrow0^{+}}N_{3}(K,\varphi_{1},\cdots,\varphi_{n})=0$$
that these two constants are what we need, and so, the claim is
proved.

\bcl\label{claim-3.8} There are constants $M_{4}(K)$ and
$N_{4}(K,\varphi_{1},\cdots,\varphi_{n})$ such that
\begin{enumerate} \item
$C_{1}(K,\varphi_{1},\cdots,\varphi_{n})\geq
M_{4}(K)-N_{4}(K,\varphi_{1},\cdots,\varphi_{n});$
\item
$\lim_{K\rightarrow1^{+}}M_{4}(K)=1$, and
\item
 $$\lim_{\|\varphi_{1}\|_{\infty}\rightarrow0^{+},\cdots,\|\varphi_{n}\|_{\infty}\rightarrow0^{+}}N_{4}(K,\varphi_{1},\cdots,\varphi_{n})=0.$$
\end{enumerate}
\ecl
By (\ref{ccww-10}), we have
$$C_{1}(K,\varphi_{1},\cdots,\varphi_{n})\geq M_{4}(K)-N_{4}(K,\varphi_{1},\cdots,\varphi_{n}),$$
where
$$M_{4}(K)=K^{-2}(Q(K))^{-2K}\frac{1}{2\pi}\int_{0}^{2\pi}|e^{it}-e^{i\theta}|^{2K-2}dt$$
and
\begin{eqnarray*}
N_{4}(K,\varphi_{1},\cdots,\varphi_{n})&=&\left(\frac{7}{6}+\frac{1}{2K^{2}}\right)\|\varphi_{1}\|_{\infty}
\\&&+\sum_{j=2}^{n}\left(\frac{47}{240}+\frac{1}{16K^{2}}\right)\|\varphi_{j}\|_{\infty}\left(\frac{3}{16}\right)^{j-2}.
\end{eqnarray*} 
The following facts
$$\lim_{K\rightarrow1^{+}}M_{4}(K)=1\;\;\mbox{and}\;\;\lim_{\|\varphi_{1}\|_{\infty}\rightarrow0^{+},\cdots,\|\varphi_{n}\|_{\infty}\rightarrow0^{+}}N_{4}(K,\varphi_{1},\cdots,\varphi_{n})=0$$
show that these two constants are what we want, and thus, the claim
is true.

Now, by the discussions of Steps \ref{bst-1} $\sim$ \ref{bst-3}, we see that the theorem is proved. \qed \medskip

\section{The proof of Proposition \ref{thm-02}}\label{sec-5}
By (\ref{eq-ch-3.0}), we have
$$f(z)=P[\varphi_{0}](z)+\sum_{k=1}^{n}(-1)^{k}G_{k}[\varphi_{k}](z),~z\in\mathbb{D}.$$
For $k\in\{1,2,\cdots,n\}$, it follows from Lemmas \ref{lem-1} and
\ref{lem-2} that $G_{k}[\varphi_{k}]$ are Lipschitz continuous in
$\mathbb{D}$. Since $P[\varphi_{0}]$ is Lipschitz continuous in
$\mathbb{D}$ if and only if the Hilbert transform of
$d\varphi_{0}(e^{i\theta})/d\theta\in L^{\infty}(\mathbb{T})$,
together with the Lipschitz continuity of $G_{k}[\varphi_{k}]$, we
conclude that $f$ is Lipschitz continuous in $\mathbb{D}$ if and
only if the Hilbert transform of
$d\varphi_{0}(e^{i\theta})/d\theta\in L^{\infty}(\mathbb{T})$, where
$k\in\{1,2,\cdots,n\}$. The proof of this proposition is complete.
\qed

\bigskip

{\bf Acknowledgements:} We are grateful to the referee for her/his comments and suggestions. This research was partly supported by the
exchange project for the third regular session of the
China-Montenegro Committee for Cooperation in Science and Technology
(No. 3-13), the Hunan Provincial Education Department Outstanding
Youth Project (No. 18B365), the Science and Technology Plan Project
of Hengyang City (No. 2018KJ125), the National Natural Science
Foundation of China (No. 11571216), the Science and Technology Plan
Project of Hunan Province (No. 2016TP1020), the Science and
Technology Plan Project of Hengyang City (No. 2017KJ183), and the
Application-Oriented Characterized Disciplines, Double First-Class
University Project of Hunan Province (Xiangjiaotong [2018]469).


\end{document}